\documentclass{svjour3}       
\journalname{}
\smartqed  % flush right qed marks, e.g. at end of proof
\usepackage{graphicx}
\usepackage{tikz}
\usetikzlibrary{patterns}
\usetikzlibrary{shapes,calc,shapes,arrows}
\usepackage{lipsum}
\usepackage{amsfonts}
\usepackage{graphicx}
\usepackage{epstopdf}
\usepackage{halloweenmath}
\usepackage{algorithm}
\usepackage{algorithmic}

\usepackage{commath,amsopn}
\usepackage{multirow}
\usepackage{stmaryrd}
\usepackage{xifthen}
\usepackage{booktabs}
\usepackage{subfigure}
\usepackage{epstopdf}
\usepackage{hyperref}
\hypersetup{colorlinks=true,linkcolor=blue,citecolor=blue}

\newcommand{\veci}{\mathbf{i}}

\newcommand{\vecr}{\mathbf{r}}

\newcommand{\vecw}{\mathbf{w}}

\newcommand{\vecmatxi}{\vec{\boldsymbol{\xi}}}
\newcommand{\vecmateta}{\vec{\boldsymbol{\eta}}}
\newcommand{\matxi}{\boldsymbol{\xi}}
\newcommand{\mateta}{\boldsymbol{\eta}}

\newcommand{\mata}{\mathbf{A}}
\newcommand{\matb}{\mathbf{B}}
\newcommand{\matc}{\mathbf{C}}

\newcommand{\matx}{\mathbf{X}}
\newcommand{\matu}{\mathbf{U}}

\newcommand{\matW}{\mathbf{W}}

\newcommand{\matA}{\mathbf{A}}
\newcommand{\matB}{\mathbf{B}}
\newcommand{\matC}{\mathbf{C}}

\newcommand{\matG}{\mathbf{G}}
\newcommand{\matI}{\mathbf{I}}
\newcommand{\matH}{\mathbf{H}}

\newcommand{\matS}{\mathbf{S}}
\newcommand{\matX}{\mathbf{X}}
\newcommand{\matU}{\mathbf{U}}

\newcommand{\matY}{\mathbf{Y}}
\newcommand{\matZ}{\mathbf{Z}}

\newcommand{\tensA}{\mathcal{A}}

\newcommand{\tensE}{\mathcal{E}}
\newcommand{\tensX}{\mathcal{X}}
\newcommand{\tensI}{\mathcal{I}}
\newcommand{\tensY}{\mathcal{Y}}
\newcommand{\tensM}{\mathcal{M}}
\newcommand{\tensN}{\mathcal{N}}
\newcommand{\tensU}{\mathcal{U}}
\newcommand{\tensS}{\mathcal{S}}
\newcommand{\tensL}{\mathcal{L}}

\newcommand{\tensH}{\mathcal{H}}

\newcommand{\subjectto}{\mathrm{s.\,t.}}
\newcommand{\T}{\mathsf{T}}

\newcommand{\vecmatW}{{\vec{\matW}}}

\newcommand{\vecmatZ}{{\vec{\matZ}}}
\newcommand{\vecmatY}{{\vec{\matY}}}

\DeclareMathOperator{\proj}{P}
\DeclareMathOperator{\tr}{tr}

\DeclareMathOperator{\grad}{grad}
\DeclareMathOperator{\tangent}{T}
\DeclareMathOperator{\rmvec}{vec}
\DeclareMathOperator{\ten}{ten}

\DeclareMathOperator*{\argmin}{arg\,min}

\begin{document}

\title{Riemannian Preconditioned Algorithms for Tensor Completion via Tensor Ring Decomposition\thanks{BG was supported by the Young Elite Scientist Sponsorship Program by CAST. YY was funded by the National Natural Science Foundation of China (grant No.12288201).}
}
% \subtitle{Do you have a subtitle?\\ If so, write it here}

\titlerunning{Riemannian Preconditioned Algorithms for Tensor Ring Completion}        % if too long for running head

\author{Bin Gao \and Renfeng Peng \and Ya-xiang Yuan}


\institute{Bin Gao \and Ya-xiang Yuan \at
    State Key Laboratory of Scientific and Engineering Computing, Academy of Mathematics and Systems Science, Chinese Academy of Sciences, Beijing, China \\
              \email{\{gaobin,yyx\}@lsec.cc.ac.cn }
           \and
           Renfeng Peng \at
           State Key Laboratory of Scientific and Engineering Computing, Academy of Mathematics and Systems Science, Chinese Academy of Sciences, and University of Chinese Academy of Sciences, Beijing, China\\
           \email{pengrenfeng@lsec.cc.ac.cn} 
}

\date{Received: date / Accepted: date}

\maketitle

\begin{abstract}
  We propose Riemannian preconditioned algorithms for the tensor completion problem via tensor ring decomposition. A new Riemannian metric is developed on the product space of the mode-2 unfolding matrices of the core tensors in tensor ring decomposition. The construction of this metric aims to approximate the Hessian of the cost function by its diagonal blocks, paving the way for various Riemannian optimization methods. Specifically, we propose the Riemannian gradient descent and Riemannian conjugate gradient algorithms. We prove that both algorithms globally converge to a stationary point. In the implementation, we exploit the tensor structure and adopt an economical procedure to avoid large matrix formulation and computation in gradients, which significantly reduces the computational cost. Numerical experiments on various synthetic and real-world datasets---movie ratings, hyperspectral images, and high-dimensional functions---suggest that the proposed algorithms have better or favorably comparable performance to other candidates.
\keywords{Tensor completion \and tensor ring decomposition \and  Riemannian optimization \and preconditioned gradient}
\PACS{15A69 \and 58C05 \and 65K05 \and 90C30}
\end{abstract}

\section{Introduction}
The tensor completion problem, as a natural generalization of the matrix completion problem, is a task of recovering a tensor based on its partially observed entries. In practice, datasets collected from real applications are often assumed to have underlying low-rank structure. Therefore, low-rank matrix decompositions are widely used in matrix completion, which can save the computational cost and storage. In the same spirit, low-rank tensor decompositions play a significant role in tensor completion; applications can be found across various fields, e.g., recommendation systems~\cite{kasai2016low,dong2022new}, image processing~\cite{liu2012tensor}, and interpolation of high-dimensional functions~\cite{steinlechner2016riemannian,khoo2021efficient}.

In this paper, we consider tensor ring decomposition and focus on the following tensor completion problem with bounded tensor ring rank $\mathbf{r}:=(r_1,\dots,r_{d})$. Given a partially observed $d$-th order tensor $\tensA\in\mathbb{R}^{n_1\times\cdots\times n_d}$ on an index set $\Omega\subseteq[n_1]\times\cdots\times[n_d]$, where $[n_k]:=\{1,2,\dots,n_k\}$ for $k=1,\dots,d$ and $d\geq 3$ is a positive interger. The tensor completion problem is formulated as follows, 
\begin{equation}\label{eq:original problem}
	\begin{array}{cc}
			\min\ &\ \frac{1}{2}\left\| \proj_\Omega(\llbracket\tensU_1,\tensU_2,\dots, \tensU_d\rrbracket)-\proj_\Omega(\tensA)\right\|_\mathrm{F}^2\\ 
		\subjectto &(\tensU_1,\tensU_2,\dots,\tensU_d)\in\tensM_\tensU, 
	\end{array}
\end{equation}
where $\llbracket\tensU_1,\tensU_2,\dots,\tensU_d\rrbracket\in\mathbb{R}^{n_1\times\cdots\times n_d}$ denotes the tensor ring decomposition with core tensors $\tensU_k\in\mathbb{R}^{r_k\times n_k\times r_{k+1}}$ for $k\in[d]$; see definitions in section~\ref{sec: preliminaries}. $\proj_\Omega$ denotes the projection operator onto $\Omega$, namely, $\proj_\Omega(\tensX)(i_1,\dots,i_d)=\tensX(i_1,\dots,i_d)$ if $(i_1,\dots,i_d)\in\Omega$, otherwise $\proj_\Omega(\tensX)(i_1,\dots,i_d)=0$. The search space of~\eqref{eq:original problem} is defined by a product space of core tensors, i.e., 
\[\tensM_\tensU:=\mathbb{R}^{r_1\times n_1\times r_2}\times \mathbb{R}^{r_2\times n_2\times r_3}\times\cdots \times \mathbb{R}^{r_d\times n_d\times r_1}.\]

\paragraph{Related works and motivation} 
Tensor completion has several different formulations, one type is based on the nuclear norm minimization. In matrix completion, the nuclear norm of a matrix---a convex relaxation of matrix rank---is minimized. Liu et al.~\cite{liu2012tensor} extended the ``nuclear norm'' to tensor by calculating the sum of nuclear norms of unfolding matrices of a tensor, and applied the alternating direction method of multipliers algorithm to solve the tensor completion problem. Since the tensor data observed in real world may be perturbed by noise, Zhao et al.~\cite{zhao2022robust} focused on the robust tensor completion problem based on the tubal nuclear norm, and proposed a proximal majorization-minimization algorithm. These algorithms require storing tensors in full size, and the number of parameters is $n_1\cdots n_d$, which scales exponentially to the dimension~$d$.

Instead of working with full-size tensors, tensor decompositions~\cite{kolda2009tensor} allow us to take advantage of the low-rank structure in tensor completion problems, and thus can reduce the number of parameters in search space and save storage. In view of the block structure in tensor decompositions, one can develop alternating minimization methods by updating one block while fixing others.
Jain and Oh~\cite{jain2014provable} proposed an alternating minimization method for symmetric third order tensors in canonical polyadic (CP) decomposition. For Tucker decomposition, the alternating minimization, also called the alternating least squares (ALS) algorithm, was considered by Andersson and Bro~\cite{andersson199893} in which the subproblem is a least squares problem. Tensor train (TT) decomposition~\cite{oseledets2010tt,oseledets2011tensor}, also known as matrix product states (MPS)~\cite{verstraete2008matrix,schollwock2011density} in computational physics, decomposes a tensor into $d$ core tensors. Grasedyck et al.~\cite{grasedyck2015variants} presented an alternating direction fitting algorithm for the tensor train completion problem. Recently, tensor ring (TR) decomposition was proposed by Zhao et al.~\cite{zhao2016tensor} as a generalization of tensor train decomposition, and is also known as MPS with periodic boundary conditions in computational physics. The ALS algorithm was developed to solve the TR decomposition problem in~\cite{zhao2016tensor}, and it was further applied to tensor completion~\cite{wang2017efficient}. In general, ALS methods in which the number of parameters scales linearly to the dimension~$d$ have been proved to be effective for the tensor completion problem. However, they may suffer from overfitting and require careful initialization in practice~\cite{chen2020tensor}.

More recently, Riemannian optimization working with tensor-based manifolds appears to be prosperous for solving completion problems. Since the search space is a manifold, one can benefit from different geometric tools and develop efficient optimization methods on the manifold where the convergence can be guaranteed; see~\cite{absil2009optimization,boumal2023intromanifolds} for an overview. Acar et al.~\cite{acar2011scalable} proposed an Euclidean nonlinear conjugate gradient method for tensor completion via CP decomposition. Dong et al.~\cite{dong2022new} proposed Riemannian gradient and Riemannian conjugate gradient methods based on a Riemannian metric on the product space of matrices in polyadic decomposition. For tensor completion in Tucker decomposition, Kressner et al.~\cite{kressner2014low} proposed a Riemannian conjugate gradient method on the manifold of rank-constrained tensors. Kasai and Mishra~\cite{kasai2016low} introduced a preconditioned metric and considered the quotient geometry of the manifold via Tucker decomposition. The corresponding Riemannian conjugate gradient algorithm was proposed. Based on geometric properties of the manifold of tensors with fixed TT rank, Steinlechner~\cite{steinlechner2016riemannian} proposed a Riemannian conjugate gradient algorithm. Furthermore, Cai et~al.~\cite{cai2022tensor} investigated the quotient geometry on this manifold, and proposed Riemannian gradient, Riemannian conjugate gradient and Riemannian Gauss--Newton algorithms. 

Tensor ring decomposition is a generalization of TT decomposition. On the one hand, it provides a flexible choice of tensor rank. Specifically, the unfolding matrices of core tensors in TR are not necessarily of full rank while the mode-1 and mode-3 unfolding matrices of core tensors in TT have to be of full rank. Moreover, TR decomposition permits a more comprehensive exploration of information along mode-$1$ and mode-$d$ by relaxing the rank constraint on the first and last cores in TT from one to an arbitrary interger, i.e., from TT rank $(1,r_2,\dots,r_d,1)$ to TR rank $(r_1,r_2,\dots,r_d)$, enabling higher compressibility and flexibility~\cite{zhao2019learning}. On the other hand, although tensors with bounded TR rank do not form a Riemannian submanifold of $\mathbb{R}^{n_1\times n_2\times \cdots\times n_d}$, TR decomposition preserves a similar block structure as TT, which allows us to endow the search space with a manifold structure and to develop a Riemannian metric that has a preconditioning effect. More precisely, we develop Riemannian optimization methods on a product space of unfolding matrices
\[\tensM:=\mathbb{R}^{n_1\times r_1r_{2}}\times\mathbb{R}^{n_2\times r_2r_{3}}\times\cdots\times\mathbb{R}^{n_{d-1}\times r_{d-1}r_{d}}\times\mathbb{R}^{n_d\times r_dr_{1}}\] 
for the tensor completion problem via TR decomposition by an equivalent reformulation of~\eqref{eq:original problem}. We refer to the search space $\tensM$ as a \emph{Riemannian manifold} by endowing $\tensM$ with a non-Euclidean metric.

\paragraph{Contributions} 
We formulate the tensor completion problem via tensor ring decomposition, in which the search space $\tensM$ is a product space of the mode-2 unfolding matrices of the core tensors in tensor ring decomposition. We design a \emph{preconditioned} metric on the search space, and propose the Riemannian gradient descent and Riemannian conjugate gradient algorithms to solve the tensor completion problem.
We prove that every accumulation point of the sequence generated by the proposed algorithms is a stationary point. To the best of our knowledge, this is the first Riemannian formulation of the tensor completion problem in TR decomposition. 

The computation of (Riemannian) gradients involves large matrix formulation and multiplication. Its computational cost is exponential to the dimension $d$, which is unaffordable in practice. In order to improve the efficiency of the proposed algorithms, we adopt an economical procedure, which has polynomial complexity to the dimension $d$, to compute the (Riemannian) gradients by exploiting the tensor structure in TR.

We compare the proposed algorithms with existing methods in various tensor completion tasks on both synthetic and real-world datasets, including movie ratings, hyperspectral images, and high-dimensional functions. The numerical results illustrate that the TR-based algorithms have better or comparable recovery performance than the others. In addition, the proposed algorithms are favorably comparable to the alternating least squares algorithm among TR-based algorithms.

\paragraph{Organization} 
First, the preliminaries of tensor ring decomposition and the tensor completion problem are introduced in section~\ref{sec: preliminaries}. Next, we develop Riemannian algorithms and present computational details in section~\ref{sec: tensor completion }. The convergence results of the proposed algorithms are shown in section~\ref{sec: convergence}. Numerical results are reported in section~\ref{sec: numerical exps}. Finally, the conclusion is given in section~\ref{sec: conclusion}.

\section{Preliminaries} \label{sec: preliminaries}
In this section, we first introduce the notation in tensor operations and the definition of tensor ring decomposition. Next, a reformulation of the TR-based tensor completion problem~\eqref{eq:original problem} is described.

The following notions are involved in tensor computations~\cite{kolda2009tensor}. First, we define an index mapping $\pi_k$ by
\begin{equation}
    \pi_k:(i_1,\dots,i_{k-1},i_{k+1},\dots,i_d)\mapsto 1 + \sum_{\ell \neq k, \ell = 1}^d(i_\ell-1)J_\ell
\end{equation}
with $J_\ell = \prod_{m = 1, m \neq k}^{\ell-1} n_m$ for $k=1,\dots,d$.
The mode-$k$ unfolding of a tensor $\mathcal{X} \in \mathbb{R}^{n_1 \times\cdots\times n_d}$ is denoted by a matrix $ \matx_{(k)}\in\mathbb{R}^{n_k\times n_{-k}} $, where $n_{-k}:=\prod_{i\neq k}n_i$. The $ (i_1,i_2,\dots,i_d)$-th entry of $\mathcal{X}$ corresponds to the $(i_k,j)$-th entry of $ \matx_{(k)} $, where
$ j = \pi_k(i_1,\dots,i_{k-1},i_{k+1},\dots,i_d)$. Similarly, the mode-$k$ unfolding of indices in an index set $\Omega$ is defined by $\Omega_{(k)}{:=}\{(i_k,\pi_k(i_1,\dots,i_{k-1},i_{k+1},\dots,i_d)):\ (i_1,\dots,i_d)\in\Omega\}$. The inner product between two tensors $\tensX,\tensY\in\mathbb{R}^{n_1\times\cdots\times n_d}$ is defined by $\langle\tensX,\tensY\rangle := \sum_{i_1=1}^{n_1} \cdots \sum_{i_d=1}^{n_d} \tensX({i_1,\dots,i_d})\tensY({i_1,\dots,i_d})$. 
The Frobenius norm of a tensor $\tensX$ is defined by $\|\tensX\|_\mathrm{F}:=\sqrt{\langle\tensX,\tensX\rangle}$.

\begin{definition}[tensor ring decomposition~\cite{zhao2016tensor}]\label{def: tensor ring}
    Given a $d$-th order tensor $\tensX\in\mathbb{R}^{n_1\times\cdots\times n_d}$, tensor ring decomposition, denoted by \[ \tensX=\llbracket\tensU_1,\dots,\tensU_d\rrbracket, \] decomposes $\tensX$ into $d$ core tensors $\tensU_k\in\mathbb{R}^{r_k\times n_k\times r_{k+1}}$ for $k=1,\dots,d$ and $r_{d+1}=r_1$. 
    Specifically, the $(i_1,\dots,i_d)$-th entry of $\tensX$ is represented by the trace of products of $d$ matrices, i.e., \[\tensX(i_1,\dots,i_d)=\tr\left(\matU_1(i_1)\cdots\matU_d(i_d)\right),\]
    where $\matu_k(i_k)=\tensU_k(:,i_k,:)\in\mathbb{R}^{r_{k}\times r_{k+1}}$ is the lateral slice matrix of $\tensU_k$ for $i_k\in[n_k]$. 
\end{definition}

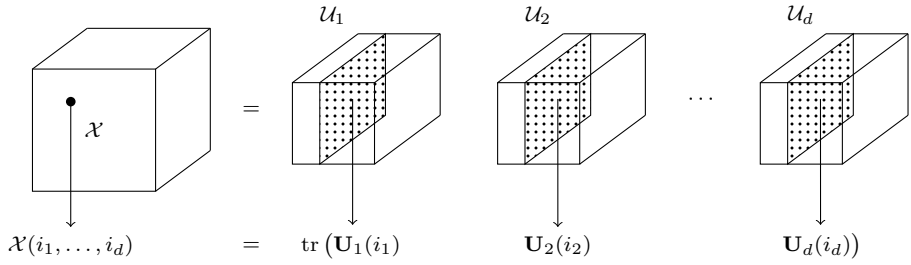
\begin{figure}[htbp]
	\centering
	\begin{tikzpicture}[scale = 1.8]
        \coordinate (base) at (-1.9,1.6); % The left vertex of the frontal square.
        \coordinate (numLevel) at (0,0.3); % Control the level of mathematical expressions
        %% Tensor X
		\draw[-] (base) rectangle ($(base)+(0.9,-0.9)$);
		\draw[-] (base) -- ($(base)+(0.4,0.3)$);
		\draw[-] ($(base)+(0.4,0.3)$) -- ($(base)+(0.4+0.9,0.3)$);
		\draw[-] ($(base)+(0.4+0.9,0.3)$) -- ($(base)+(0.9,0)$);
		\draw[-] ($(base)+(0.4+0.9,0.3)$) -- ($(base)+(0.4+0.9,-0.6)$);
		\draw[-] ($(base)+(0.4+0.9,0.3-0.9)$) -- ($(base)+(0.9,-0.9)$);
		\node (v) at ($(base)+(0.45,-0.45)$) {$\tensX$};

        \fill ($(base)+(0.28,-0.24)$) circle (1pt);
		\node (x123) at ($(base)+(0.28,-1.6)+(numLevel)$) {$\tensX(i_1,\dots,i_d)$};
		\draw[->] ($(base)+(0.28,-0.24)$) -- (x123);
	
		\node (u) at (-0.3,1.3) {$=$};
		\node (u) at (-0.3,0.3) {$=$};

        %% Tensor U1
        \coordinate (U1) at (0,1.5); 
        \coordinate (squareSize) at (0.6,-0.6);
		\draw[-] (1.5+0.6-1.5,1.5) -- (1.98+0.6-1.5,1.86);
		\draw[-] (1.98+0.6-1.5,1.86) -- (1.38+0.6-1.5,1.86);
		\draw[-] (0.9+0.6-1.5,1.5) -- (1.38+0.6-1.5,1.86);
		\draw[-] (1.98+0.6-1.5,1.86) -- (1.98+0.6-1.5,1.26);
		\draw[-] (1.5+0.6-1.5,0.9) -- (1.98+0.6-1.5,1.26);
		\draw[-] (U1) rectangle ($(U1)+(squareSize)$);
		\node (G) at (1.2+0.6-1.5,2) {$\tensU_1$};
		\path[draw, pattern=dots] (0.2+1.5-1.5,0.9) -- (0.68+1.5-1.5,1.26) -- (0.68+1.5-1.5,1.86) -- (0.2+1.5-1.5,1.5) -- cycle;
		\node (G1) at (0.44+1.5-1.5,0.3) {$\tr\left(\matU_1(i_1)\right.$};
		\draw[->] (0.44+1.5-1.5,1.38) -- (G1);

		\draw[-] (1.5+0.6,1.5) -- (1.98+0.6,1.86);
		\draw[-] (1.98+0.6,1.86) -- (1.38+0.6,1.86);
		\draw[-] (0.9+0.6,1.5) -- (1.38+0.6,1.86);
		\draw[-] (1.98+0.6,1.86) -- (1.98+0.6,1.26);
		\draw[-] (1.5+0.6,0.9) -- (1.98+0.6,1.26);
		\draw[-] (0.9+0.6,1.5) rectangle (1.5+0.6,0.9);
		\node (G) at (1.2+0.6,2) {$\tensU_2$};
		\path[draw, pattern=dots] (0.2+1.5,0.9) -- (0.68+1.5,1.26) -- (0.68+1.5,1.86) -- (0.2+1.5,1.5) -- cycle;
		\node (G2) at (0.44+1.5,0.3) {$\matU_2(i_2)$};
		\draw[->] (0.44+1.5,1.38) -- (G2);
		
		\node(cdots) at (3,1.38) {$\cdots$};

		\draw[-] (1.5+0.6+1.92,1.5) -- (1.98+0.6+1.92,1.86);
		\draw[-] (1.98+0.6+1.92,1.86) -- (1.38+0.6+1.92,1.86);
		\draw[-] (0.9+0.6+1.92,1.5) -- (1.38+0.6+1.92,1.86);
		\draw[-] (1.98+0.6+1.92,1.86) -- (1.98+0.6+1.92,1.26);
		\draw[-] (1.5+0.6+1.92,0.9) -- (1.98+0.6+1.92,1.26);
		\draw[-] (0.9+0.6+1.92,1.5) rectangle (1.5+0.6+1.92,0.9);
		\node (G) at (1.2+0.6+1.92,2) {$\tensU_d$};
		\path[draw, pattern=dots] (0.2+1.5+1.92,0.9) -- (0.68+1.5+1.92,1.26) -- (0.68+1.5+1.92,1.86) -- (0.2+1.5+1.92,1.5) -- cycle;
		\node (Gd1) at (0.44+1.5+1.92,0.3) {$\left.\matU_d(i_{d})\right)$};
		\draw[->] (0.44+1.5+1.92,1.38) -- (Gd1);
	\end{tikzpicture}
	\caption{Illustration of tensor ring decomposition of a tensor }
	\label{fig: Tensor Ring}
\end{figure}

The tuple $\mathbf{r}=(r_1,\dots,r_d)$ is referred to as the \emph{tensor ring rank} (TR rank). Figure~\ref{fig: Tensor Ring} illustrates the TR decomposition of a tensor.

In view of the cyclic symmetry of trace operator, we rewrite the $(i_1,\dots,i_d)$-th entry of tensor $\tensX$ in tensor ring decomposition as follows, 
\begin{equation*}
    \begin{aligned}
        &\tr(\matU_1(i_1)\cdots\matU_d(i_d))=\tr\left(\matU_k(i_k)\cdots\matU_d(i_d)\matU_1(i_1)\cdots\matU_{k-1}(i_{k-1})\right)\nonumber\\
        &=\langle\rmvec(\matU_k(i_k)),\rmvec((\matU_{k+1}(i_{k+1})\cdots \matU_{d}(i_{d})\matU_{1}(i_{1})\cdots \matU_{k-1}(i_{k-1}))^\T)\rangle,
    \end{aligned}
\end{equation*}
where $\rmvec(\cdot)$ denotes the column vectorization of a matrix. In fact, $\rmvec(\matU_k(i_k))^\T$ is the $i_k$-th row of the mode-$2$ unfolding matrix $\left(\tensU_k\right)_{(2)}\in\mathbb{R}^{n_k\times r_kr_{k+1}}$ of the core tensor $\tensU_k$. Additionally, given a matrix $\matW_k\in\mathbb{R}^{n_k\times (r_kr_{k+1})}$ for fixed $n_1,\dots,n_d$ and $r_1,\dots,r_d$, the second tensorization operator maps  $\matW_k$ to a tensor $\ten_{(2)}(\matW_k)\in\mathbb{R}^{r_k\times n_k\times r_{k+1}}$ defined by 
\[\ten_{(2)}(\matW_k)(i_1,i_2,i_3):=\matW_k(i_2,i_1+(i_3-1)r_k)\]
for $(i_1,i_2,i_3)\in[r_k]\times[n_k]\times[r_{k+1}]$. We observe that $(\ten_{(2)}(\matW_k))_{(2)}=\matW_k$ holds. Therefore, the second tensorization operator is invertible. The definition of subchain tensors is given as follows. 

\begin{definition}[subchain tensor]\label{def: subchain}
    The subchain tensor $\tensU_{\neq k}\in\mathbb{R}^{r_{k}\times n_{-k}\times r_{k+1}}$ is defined by its lateral slice matrices, i.e.,
    \begin{equation}
        \label{eq: def U neq k} \matU_{\neq k}\left(\pi_k(i_1,\dots,i_{k-1},i_{k+1},\dots,i_d)\right){:=}\left(\prod_{j=k+1}^{d}\matU_{j}(i_{j})\prod_{j=1}^{k-1}\matU_{j}(i_{j})\right)^\T
        %\left(\matU_{k+1}(i_{k+1})\cdots \matU_{d}(i_{d})\matU_{1}(i_{1})\cdots \matU_{k-1}(i_{k-1})\right)^\T.
    \end{equation}
    for $(i_1,\dots,i_{k-1},i_{k+1},\dots,i_d)\in[n_1]\times\cdots\times[n_{k-1}]\times[n_{k+1}]\times\cdots\times[n_d]$ and $k\in[d]$.
\end{definition}

Given the mode-2 unfolding matrix $\left(\tensU_k\right)_{(2)}$ and subchain tensor in Definition~\ref{def: subchain}, the mode-$k$ unfolding of a tensor $\tensX$ in tensor ring decomposition equals to the product of two smaller matrices (see~\cite[Theorem 3.5]{zhao2016tensor}), i.e., 
\[\matX_{(k)}=\matW_k^{}\matW_{\neq k}^\T,\] 
where  $\matW_k:=\left(\tensU_k\right)_{(2)}\in\mathbb{R}^{n_k\times r_kr_{k+1}}$ and $\matW_{\neq k}:=\left(\tensU_{\neq k}\right)_{(2)}\in\mathbb{R}^{n_{-k}\times r_kr_{k+1}}$ for $k=1,\dots,d$. In the light of this fact, we are able to formulate the TR-based tensor problems on {mode-2} unfolding matrices $\matW_1,\dots,\matW_d$. Specifically, the update of a TR tensor $\tensX=\llbracket\tensU_1,\dots,\tensU_d\rrbracket$ involves the updates of each core tensors $\tensU_1,\dots,\tensU_d$, which is equivalent to updating the mode-2 unfolding matrices~$\matW_1,\dots,\matW_d$ by matricizations $\matW_k=(\tensU_k)_{(2)}$. In other words, mode-$2$ unfolding matrices are adequate for all tensor-related computations. It is worth noting that given $\matW_1,\dots,\matW_d$, the matrices $\matW_k^{}\matW_{\neq k}^\T$ and $\matW_j^{}\matW_{\neq j}^\T$ represent a same tensor after their respective tensorizations for $j\neq k$. However, we never compute $\matW_k^{}\matW_{\neq k}^\T$ explicitly in practice.

Subsequently, the tensor completion problem~\eqref{eq:original problem} can be reformulated on the product space of the mode-2 unfolding matrices of the core tensors in TR decomposition \[\tensM=\mathbb{R}^{n_1\times r_1r_{2}}\times\mathbb{R}^{n_2\times r_2r_{3}}\times\cdots\times\mathbb{R}^{n_{d-1}\times r_{d-1}r_{d}}\times\mathbb{R}^{n_d\times r_dr_{1}}.\]
The Frobenius norm on $\tensM$ is defined by $\|\vecmatW\|_\mathrm{F}:=\sqrt{\sum_{k=1}^d\|\matW_k\|_\mathrm{F}^2}$ for $\vecmatW=\left(\matW_1,\dots,\matW_d\right)\in\tensM$. 

In this paper, we focus on the following tensor completion problem combined with a sampling set $\Omega$ and a regularization term $r$,
\begin{align} %
    \min_{\vec{\matW}=\left(\matW_1,\dots,\matW_d\right)\in\tensM}&\ f(\vecmatW):=f_{\Omega}(\vecmatW)+r(\vecmatW).
    \label{eq: newly stated problem}
\end{align}
The objective function $f$ consists of two parts: one is the cost function 
\begin{equation*}
    f_{\Omega}(\vecmatW):=\frac{1}{2p}\left\| \proj_\Omega(\tensX)-\proj_\Omega(\tensA)\right\|_\mathrm{F}^2=\frac{1}{2p}\left\|\proj_{\Omega_{(k)}} \left(\matW_k^{}\matW_{\neq k}^\T-\matA_{(k)}\right) \right\|_\mathrm{F}^2,
    % \frac{1}{2p}\left\|\proj_{\Omega_{(k)}} \left(\matW_k^{}\matW_{\neq k}^\T\right) - \proj_{\Omega_{(k)}}\left(\matA_{(k)}\right) \right\|_\mathrm{F}^2,
\end{equation*}
where $p:=\frac{\lvert \Omega \rvert}{n_1\cdots n_d}$ 
is the \emph{sampling rate}, $\mata_{(k)}$ denotes the mode-$k$ unfolding matrix of $\tensA$, and 
$\Omega_{(k)}$ is the mode-$k$ unfolding of the sampling set $\Omega$; 
the other is a regularization term $r(\vecmatW)$ and we choose
$$r(\vecmatW):={\frac{\lambda}{2}\|\vecmatW\|_\mathrm{F}^2}$$ in a similar fashion as Maximum-margin Matrix Factorization~\cite{srebro2004maximum} with $\lambda>0$. The regularization term keeps the variable $\vecmatW$ in a compact subset of $\tensM$ and guarantees the convergence; see analysis in section~\ref{sec: convergence}.

Since the second tensorization operator $\ten_{(2)}(\cdot)$ is invertible when $n_1,\dots,n_d$ and $r_1,\dots,r_d$ are fixed, it turns out that the search space in~\eqref{eq: newly stated problem} is equivalent to the one in the original problem~\eqref{eq:original problem} and the two search spaces are connected by matricizations $\matW_k=(\tensU_k)_{(2)}$ and tensorizations $\tensU_k=\ten_{(2)}(\matW_k)$ for $(\matW_1,\dots,\matW_d)\in\tensM$ and $(\tensU_1,\dots,\tensU_d)\in\tensM_\tensU$, i.e., 
\[
    \begin{array}{c}
        (\tensU_1,\dots,\tensU_d)\in\tensM_\tensU \\[3pt]
        \mathbb{R}^{r_1\times n_1\times r_2}\times\cdots \times \mathbb{R}^{r_d\times n_d\times r_1}
    \end{array}
    \ 
    \begin{array}{c}
        \text{matricization} \\[2pt]
        \underscriptrightarrow{\ \ \ \ \ \ \matW_k=(\tensU_k)_{(2)}\ \ \ \ \ \ }\\[5pt]
        \overscriptleftarrow{\ \ \ \ \tensU_k=\ten_{(2)}(\matW_k)\ \ \ \ }\\[2pt]
        \text{tensorization}
    \end{array}
    \ 
    \begin{array}{c}
        \tensM\ni(\matW_1,\dots,\matW_d) \\[3pt]
        \mathbb{R}^{n_1\times r_1r_{2}}\times\cdots\times\mathbb{R}^{n_d\times r_dr_{1}}
    \end{array}.
\]
Additionally, we observe that an element $\vec{\matW}=\left(\matW_1,\dots,\matW_d\right)$ in $\tensM$ can also represent an element $\llbracket\ten_{(2)}(\matW_1),\dots,\ten_{(2)}(\matW_d)\rrbracket\in\mathbb{R}^{n_1\times\cdots\times n_d}$ by TR decomposition. We define the mapping 
\begin{equation}
    \tau:\tensM\to\mathbb{R}^{n_1\times\cdots\times n_d}\label{eq: mapping tau},\ \tau(\vecmatW){:=}\llbracket\ten_{(2)}(\matW_1),\dots,\ten_{(2)}(\matW_d)\rrbracket,
\end{equation}
which is adopted to generate systhetic data in section~\ref{sec: numerical exps}.

Given the objective function $f$, the first-order derivative with respect to $\matW_k$ has the following form,  
\begin{equation*}
    \partial_{\matW_k} f(\vecmatW)=\partial_{\matW_k} f_\Omega(\vecmatW)+\partial_{\matW_k} r(\vecmatW)=\frac{1}{p} \matS_{(k)}\matW_{\neq k}+\lambda \matW_{k}\quad\text{for\ } k=1,\dots,d,
\end{equation*}
where $\tensS:=\proj_\Omega(\tau(\vecmatW))-\proj_\Omega(\tensA)$ is called the residual tensor and $\matS_{(k)}$ is the mode-$k$ unfolding matrix of tensor $\tensS$. Therefore, the Euclidean gradient of $f_\Omega$ at $\vecmatW\in\tensM$ can be computed as follows, 
\begin{equation}
    \nabla f_{\Omega}(\vecmatW)=\left(\frac{1}{p} \matS_{(1)}\matW_{\neq 1}, \frac{1}{p} \matS_{(2)}\matW_{\neq 2}, \dots, \frac{1}{p} \matS_{(d)}\matW_{\neq d}\right).
    \label{eq: Euclidean gradient of f Omega} 
\end{equation}
Moreover, one can develop Euclidean gradient descent algorithms (e.g.,~\cite{yuan2018higher}) to solve the completion problem~\eqref{eq: newly stated problem}. 
Recently, Riemannian preconditioned algorithms are considered in which the search space is endowed with a non-Euclidean metric. The construction of this metric aims to approximate the Hessian of the cost function by its ``diagonal blocks''. These algorithms improve the performance of Euclidean-based algorithms, and are successfully applied to matrix and tensor completion (e.g.,~\cite{mishra2012riemannian,kasai2016low,dong2022new,cai2022tensor}). However, the extension to TR is not straightforward since it involves large matrix formulation and computation. In next section, we consider how to develop efficient preconditioned algorithms for the tensor completion problem~\eqref{eq: newly stated problem}.

\section{Tensor completion algorithms}\label{sec: tensor completion }
We first develop a preconditioned metric on the manifold $\tensM$. The corresponding Riemannian gradient is derived under this metric. Next, we propose the Riemannian gradient descent and Riemannian conjugate gradient algorithms. An efficient procedure for computing the Riemannian gradient is proposed in the end.

\subsection{A preconditioned metric}\label{subsec: preconditioned metric}
The idea of developing a preconditioned metric on $\tensM$ is to take advantage of the second-order information of the cost function $f_\Omega$ and to formulate a search direction that approximates the Newton direction. Specifically, we intend to construct an operator $\tensH(\vecmatW):\tangent_\vecmatW\tensM\to \tangent_\vecmatW\tensM$ such that 
\begin{equation}
    \langle\tensH(\vecmatW)[\vecmatxi],\vecmateta\rangle\approx\nabla^2 f_{\Omega}(\vecmatW)[\vecmatxi,\vecmateta]\label{eq: approx}
\end{equation}
for all $\vecmatxi,\vecmateta\in \tangent_\vecmatW\tensM\simeq\tensM$, 
where $\tangent_\vecmatW\tensM$ denotes the tangent space to $\tensM$ at $\vecmatW
\in\tensM$ and $\nabla^2 f_{\Omega}$ denotes the (Euclidean) Hessian of $f_\Omega$. Note that 
    \begin{equation*}
        \begin{aligned}
            \tangent_\vecmatW\tensM&=\tangent_{\matW_1}\mathbb{R}^{n_1\times r_1r_{2}}\times\cdots\times\tangent_{\matW_{d-1}}\mathbb{R}^{n_{d-1}\times r_{d-1}r_{d}}\times\tangent_{\matW_d}\mathbb{R}^{n_d\times r_dr_{1}}\\
            &=\mathbb{R}^{n_1\times r_1r_{2}}\times\cdots\times\mathbb{R}^{n_{d-1}\times r_{d-1}r_{d}}\times\mathbb{R}^{n_d\times r_dr_{1}}.
        \end{aligned}
    \end{equation*}
    Therefore, a tangent vector $\vecmatxi\in\tangent_\vecmatW\tensM$ can be expressed by $\vecmatxi=(\matxi_1,\matxi_2,\dots,\matxi_d)$ with $\matxi_k\in\mathbb{R}^{n_k\times r_kr_{k+1}}$.

To this end, we start from computing the explicit form of the Hessian operator
\[\nabla^2 f_\Omega(\vecmatW)[\vecmatxi,\vecmateta]=\sum_{k=1}^{d}\langle\partial_{\matW_k,\matW_k}^2f_\Omega(\vecmatW)[\vecmatxi],\vecmateta\rangle+\sum_{\ell,m=1,\ell\neq m}^d\langle\partial_{\matW_\ell,\matW_m}^2f_\Omega(\vecmatW)[\vecmatxi],\vecmateta\rangle\] for all $\vecmatxi,\vecmateta\in \tangent_\vecmatW\tensM$. By direct calculations, the ``diagonal blocks'' of $\nabla^2 f_\Omega(\vecmatW)$ have the following forms,
\begin{equation}\label{eq: 2nd order derivative}
    \partial_{\matW_k,\matW_k}^2 f_{\Omega}(\vecmatW)[\vecmatxi]=\frac{1}{p} \proj_{\Omega_{(k)}}\left(\matxi_k\matW_{\neq k}^\T \right)\matW_{\neq k}^{}\quad\text{for } k=1,\dots,d.
\end{equation} 
Since computing the ``off-diagonal blocks'' $\partial_{\matW_\ell,\matW_m}^2f_\Omega(\vecmatW)$ is complicated, we consider only the ``diagonal blocks'' as a trade-off between accuracy and computational cost to form an operator $\tensH(\vecmatW)$. An instinctive approach to construct the operator is to directly apply the second-order derivatives, i.e., 
\[\tensH_{\Omega}(\vecmatW)[\vecmatxi]:=\left(\partial_{\matW_1,\matW_1}^2 f_{\Omega}(\vecmatW)[\vecmatxi],\dots,\partial_{\matW_d,\matW_d}^2 f_{\Omega}(\vecmatW)[\vecmatxi]\right),\]
which relies on a specific sampling set $\Omega$. We intend to design an operator that is applicable to a class of tensor completion problem with a subsampling pattern by taking expectation on~\eqref{eq: 2nd order derivative}.
More accurately, we suppose that the indices in $\Omega$ are i.i.d. samples from the Bernoulli distribution with probability $p$. We eliminate the projection operator $\proj_\Omega$ by taking expectation on ``diagonal blocks'' over~$\Omega$. Hence, we can define the operator $\tensH(\vecmatW)$ by 
\begin{equation}
    \label{eq: the explicit representation of H(U)}
    \begin{aligned}
        \tensH(\vecmatW)[\vecmatxi]&:=\left(\mathbb{E}_\Omega\left[\partial_{\matW_1,\matW_1}^2 f_{\Omega}(\vecmatW)[\vecmatxi]\right],\dots,\mathbb{E}_\Omega\left[\partial_{\matW_d,\matW_d}^2 f_{\Omega}(\vecmatW)[\vecmatxi]\right]\right)\\
        &\:=\left(\matxi_1\matW_{\neq 1}^\T \matW_{\neq 1}^{}, \dots, \matxi_d\matW_{\neq d}^\T \matW_{\neq d}^{}
        \right).
    \end{aligned}
\end{equation}
Note that both $\tensH_{\Omega}$ and $\tensH$ approximate the second-order information of the cost function. The operator $\tensH$ is an expectation form of $\tensH_{\Omega}$. We adopt the operator~\eqref{eq: the explicit representation of H(U)} for computational convenience but one can consider other subsampling patterns in practice. Consequently, we define a new metric as follows.

\begin{definition}[preconditioned metric]
    $g$ is an inner product on $\tensM$ defined by
    \begin{equation}
        \label{eq: a new metric}
        g_\vecmatW(\vecmatxi,\vecmateta):=\sum_{k=1}^d\tr\left(\matxi_k\matH_k(\vecmatW)\left(\mateta_k\right)^\T\right)
    \end{equation} 
    for all $\vecmatW\in\tensM$ and $\vecmatxi,\vecmateta\in \tangent_\vecmatW \tensM$,  
    where $\matH_k(\vecmatW)\in\mathbb{R}^{r_kr_{k+1}\times r_kr_{k+1}}$ is a matrix defined by \[\matH_k(\vecmatW):=\matW_{\neq k}^\T \matW_{\neq k}^{}+\delta\matI_{r_kr_{k+1}}\] with a constant parameter $\delta>0$ and the identity matrix $\matI_{r_kr_{k+1}}\in\mathbb{R}^{r_kr_{k+1}\times r_kr_{k+1}}$.
\end{definition}

Since the matrix $\matW_{\neq k}^\T \matW_{\neq k}^{}$ is not necessarily positive definite, a shifting term $\delta\matI_{r_kr_{k+1}}$ is added to avoid singularity. Moreover, $g$ is smooth on $\tensM$, and thus is a well-defined Riemannian metric on $\tensM$. It turns out that $\tensM$ is a Riemannian manifold endowed with $g$ and the norm of a tangent vector $\vecmatxi\in \tangent_\vecmatW\tensM$ can be defined by $\|\vecmatxi\|_\vecmatW:=\sqrt{g_\vecmatW(\vecmatxi,\vecmatxi)}$. The Riemannian gradient~\cite[Sect. 3.6]{absil2009optimization}, $\grad\!f(\vecmatW)$, is the unique element in $\tangent_\vecmatW\tensM$ that satisfies 
\[{g}_\vecmatW\left(\grad\!f(\vecmatW),\vecmatxi\right)=\mathrm{D}f(\vecmatW)[\vecmatxi]:=\langle\nabla f(\vecmatW),\vecmatxi\rangle\]
for all $\vecmatxi\in \tangent_\vecmatW\tensM$. Note that 
\[g_\vecmatW(\vecmateta,\vecmatxi)=\sum_{k=1}^d\tr\left(\mateta_k\matH_k(\vecmatW)\left(\matxi_k\right)^\T\right)=\langle(\tensH+\delta\tensI)(\vecmatW)[\vecmateta],\vecmatxi\rangle,\]
where $\tensI(\vecmatW):\tangent_\vecmatW\tensM\to \tangent_\vecmatW\tensM$ is the identity operator. It turns out that \[\grad\!f(\vecmatW)=(\tensH+\delta\tensI)^{-1}(\vecmatW)[\nabla f(\vecmatW)].\] 
In view of~\eqref{eq: newly stated problem} and~\eqref{eq: approx}, $\grad\!f$ is an approximation of the Newton direction of $f$. As a result, the metric~\eqref{eq: a new metric} has a preconditioning effect on the Euclidean gradient. Therefore, we refer to the new metric~\eqref{eq: a new metric} as a \emph{preconditioned metric} on~$\tensM$. In summary, the Riemannian gradient can be computed from~\eqref{eq: Euclidean gradient of f Omega} as follows. 
\begin{proposition}[Riemannian gradient]
    The Riemannian gradient of $f$ at $\vecmatW\in\tensM$ with respect to the metric $g$ is 
    \begin{equation}\label{eq: Riemannian Gradient}
        \grad\!f(\vecmatW) = \left(\partial_{\matW_1} f(\vecmatW)\matH_1^{-1}(\vecmatW), \dots, \partial_{\matW_d} f(\vecmatW)\matH_d^{-1}(\vecmatW)\right).
    \end{equation}
\end{proposition}

\subsection{Riemannian preconditioned algorithms} 
Using the preconditioned metric~\eqref{eq: a new metric} and the Riemannian gradient~\eqref{eq: Riemannian Gradient}, we propose the Riemannian gradient descent and Riemannian conjugate gradient algorithms to solve the tensor completion problem~\eqref{eq: newly stated problem}. 

\paragraph{Riemannian gradient descent}
The Riemannian gradient descent algorithm is listed in Algorithm~\ref{alg: TR-RGD}. Note that the retraction on the manifold $\tensM$ is the identity map. For the selection of stepsize, we consider the following two strategies: 1) we adopt exact line search by solving the following optimization problem
\begin{equation}
    s^{(t)}_\mathrm{exact}:=\argmin\limits_{s>0}\, h(s)=f(\vecmatW^{(t)}+s\vecmateta^{(t)}).\label{eq: linmin}
\end{equation} 
Since $h$ is a polynomial of $s$ with degree $2d$, the solution $s^{(t)}_\mathrm{exact}$ is the root of the polynomial $h^\prime(s)$ with $2d-1$ degree;  
2) alternatively, we consider Armijo backtracking line search. Given the initial stepsize $s_0^{(t)}>0$, find the smallest integer $\ell$, such that for {$s^{(t)}=\rho^\ell 
s_0^{(t)}>{s}_{\min}$}, the inequality 
\begin{equation}
    f(\vecmatW^{(t)})-f(\vecmatW^{(t)}+s^{(t)}\vecmateta^{(t)})\geq -s^{(t)} a g_{\vecmatW^{(t)}}\left(\grad\!f(\vecmatW^{(t)}),\vecmateta^{(t)}\right)\label{eq: Armijo}
\end{equation}
holds, where $\rho, a\in(0,1), {{s}_{\min}}>0$ are backtracking parameters. The Riemannian Barzilai--Borwein (RBB) stepsize~\cite{iannazzo2018riemannian}, defined by
\begin{equation}
    \label{eq: alternating BB stepsizes}
    s^{(t)}_\mathrm{RBB1}:=\frac{\|\vecmatZ^{(t-1)}\|_{\vecmatW^{(t)}}^2}{\lvert g_{\vecmatW^{(t)}}(\vecmatZ^{(t-1)},\vecmatY^{(t-1)})\rvert}\quad\text{or}\quad s^{(t)}_\mathrm{RBB2}:=\frac{\lvert g_{\vecmatW^{(t)}}(\vecmatZ^{(t-1)},\vecmatY^{(t-1)})\rvert}{\|\vecmatY^{(t-1)}\|_{\vecmatW^{(t)}}^2},
\end{equation}
where $\vecmatZ^{(t-1)}:=\vecmatW^{(t)}-\vecmatW^{(t-1)}$ and $\vecmatY^{(t-1)}:=\grad\!f(\vecmatW^{(t)})-\grad\!f(\vecmatW^{(t-1)})$, appears to be favorable in many applications. Therefore, we set RBB to be the initial stepsize $s_0^{(t)}$. It is worth noting that the initial stepsize $s_0^{(0)}$ can also be generated by exact line search~\eqref{eq: linmin} in practice.

\begin{algorithm}[H]
    \caption{Riemannian Gradient Descent Algorithm (TR-RGD)}
    \begin{algorithmic}[1]
        \REQUIRE $f:\tensM\to\mathbb{R}$, $\vecmatW^{(0)}\in\tensM$, tolerance $\varepsilon>0$, $t=0$; {backtracking parameters $\rho, a\in(0,1), {{s}_{\min}}>0$}.
        \WHILE{the stopping criteria are not satisfied}
            \STATE Compute search direction $\vecmateta^{(t)}= -\grad\!f(\vecmatW^{(t)})$.
            \STATE Compute stepsize $s^{(t)}$ by exact line search~\eqref{eq: linmin} or Armijo backtracking~\eqref{eq: Armijo}.
            \STATE Update $\vecmatW^{(t+1)}= \vecmatW^{(t)}+s^{(t)}\vecmateta^{(t)};\ t= t+1$.
        \ENDWHILE
        \ENSURE $\vecmatW^{(t)}\in\tensM$.
    \end{algorithmic}\label{alg: TR-RGD}
\end{algorithm}

We illustrate the connection and difference between the alternating least squares algorithm for tensor ring completion~\cite{wang2017efficient} (TR-ALS) and the proposed TR-RGD algorithm in the following remark.
\begin{remark}
    Both TR-ALS and TR-RGD can be viewed as line search methods. In contrast with the TR-ALS, where $\matW_1,\dots,\matW_d$ are updated sequentially with exact line search, $(\matW_1,\dots,\matW_d)$ are treated as one point $\vecmatW\in\tensM$ in the proposed TR-RGD and updated by the Riemannian gradient descent algorithm. The proposed TR-RGD benefits from a preconditioned metric~\eqref{eq: a new metric} which is exquisitely tailored for the tensor completion problem~\eqref{eq: newly stated problem}. Moreover, TR-RGD allows a more flexible choice of stepsize rules, e.g., exact line search and Riemannian BB stepsize. Therefore, the TR-RGD method can be potentially competitive in practice.
\end{remark}

\paragraph{Riemannian conjugate gradient}
The Riemannian conjugate gradient algorithm is given in Algorithm~\ref{alg: TR-RCG}. 
For CG parameter $\beta^{(t)}$, 
we consider the Riemannian version~\cite{boumal2014manopt} of the modified Hestenes--Stiefel rule (HS+)~\cite{hestenes1952methods}  
\begin{equation}\label{eq: HS+}
    \beta^{(t)}{:=}\max\left\{\frac{g_{\vecmatW^{(t)}}\left(\grad\!f(\vecmatW^{(t)})-\grad\!f(\vecmatW^{(t-1)}),\grad\!f(\vecmatW^{(t)})\right)}{g_{\vecmatW^{(t)}}\left(\grad\!f(\vecmatW^{(t)})-\grad\!f(\vecmatW^{(t-1)}),\vecmateta^{(t-1)}\right)},0\right\}.
\end{equation}
The stepsize in Algorithm~\ref{alg: TR-RCG} is determined by Armijo backtracking line search~\eqref{eq: Armijo}. 

\begin{algorithm}[H]
    \caption{Riemannian Conjugate Gradient Algorithm (TR-RCG)}\label{alg: TR-RCG}
    \begin{algorithmic}[1]
        \REQUIRE $f:\tensM\to\mathbb{R}$, $\vecmatW^{(0)}\in\tensM$, tolerance $\varepsilon>0$, $t=0$, $\vecmateta^{(-1)}=0$; {backtracking parameters $\rho, a\in(0,1), {{s}_{\min}}>0$}.
        \WHILE{the stopping criteria are not satisfied}
            \STATE Compute search direction $\vecmateta^{(t)}= -\grad\!f(\vecmatW^{(t)})+\beta^{(t)}\vecmateta^{(t-1)}$ with CG parameter~\eqref{eq: HS+}. 
            \STATE Compute stepsize $s^{(t)}$ by Armijo backtracking line search~\eqref{eq: Armijo}.
            \STATE Update $\vecmatW^{(t+1)}= \vecmatW^{(t)}+s^{(t)}\vecmateta^{(t)};\ t= t+1$.
        \ENDWHILE
        \ENSURE $\vecmatW^{(t)}\in\tensM$.
    \end{algorithmic}
\end{algorithm}

In both algorithms, the computation of Riemannian gradient, involving large matrix formulation and multiplication, dominates the total cost. Since this cost is exponential to $d$, a~straightforward implementation is not affordable in practice. To this end, we have to come up with a procedure that can compute the Riemannian gradient efficiently. 

\subsection{Efficient computation for gradients}\label{subsec: computational skills}
We investigate the computational details of the Riemannian gradient~\eqref{eq: Riemannian Gradient} in this subsection. Generally, the computation of $\grad\!f(\vecmatW)=(\mateta_1,\dots,\mateta_d)$ involves two steps: the first is computing the Euclidean gradient $\nabla f(\vecmatW)=(\matG_1,\dots,\matG_d)$ in~\eqref{eq: Euclidean gradient of f Omega}; the second is assembling the Riemannian gradient $\grad\!f(\vecmatW)$, where
\begin{equation*}
    \begin{aligned}
        \matG_k&:=\partial_{\matW_k}f(\vecmatW)=\frac{1}{p}\matS_{(k)}\matW_{\neq k}+\lambda\matW_k,\\
        \mateta_k&:=\matG_k\matH_k^{-1}(\vecmatW)=\matG_k(\matW_{\neq k}^\T \matW_{\neq k}^{}+\delta\matI_{r_kr_{k+1}}^{})^{-1}
    \end{aligned}
\end{equation*}
for $k=1,\dots,d$. 

There are five operations in straightforward calculating the gradients:
1) form the matrix $\matW_{\neq k}$ for $k=1,\dots,d$, which requires $2\sum_{k=1}^dn_{-k}R_k$ flops, where $R_k:=r_k\left(\sum_{j=1,j\neq k}^dr_{j}r_{j+1}\right)$; 2) compute the sparse tensor $\tensS$, requiring $2\lvert \Omega \rvert r_1r_2$ flops; 3) compute the Euclidean gradient by sparse-dense matrix product which involves $2\lvert \Omega \rvert\bar{r}$ flops, where $\bar{r}:=\sum_{k=1}^d r_kr_{k+1}$; 4) compute $\matH_k(\vecmatW)$ by a dense-dense matrix multiplication with $2\sum_{k=1}^d n_{-k}(r_kr_{k+1})^2$ flops; 5) compute the Riemannian gradient through Cholesky decomposition for linear systems, which requires $\sum_{k=1}^d \left(2n_k(r_kr_{k+1})^2+C_\mathrm{chol}(r_kr_{k+1})^3\right)$ flops.
Following above operations, the straightforward computation of the Riemannian gradient totally requires 
\[2\lvert \Omega \rvert d\bar{r}+2\lvert \Omega \rvert r_1r_2+\sum_{k=1}^d \left(2n_{-k}(R_k+(r_kr_{k+1})^2)+2n_k(r_kr_{k+1})^2+C_\mathrm{chol}(r_kr_{k+1})^3\right)\]
flops. If $n_1=\cdots=n_d=n$ and $r_1=\cdots=r_d=r$, it boils down to \[2(d+1)\lvert \Omega \rvert r^2+2d(d-1)n^{d-1}r^3+2dn^{d-1}r^4+2dnr^4+C_{\mathrm{chol}}dr^6\] flops in total. In practice, the terms with order of $\mathcal{O}(n^{d-1})$ dominate the computational cost.

By using the Kronecker product structure of $\matW_{\neq k}^\T \matW_{\neq k}^{}$, the total cost can be significantly reduced. Algorithm~\ref{alg: compute uneqk t uneqk} illustrates how to compute $\matW_{\neq k}^\T \matW_{\neq k}^{}$  without forming $\matW_{\neq k}$ explicitly. The matrix multiplication is
\begin{equation}
    \label{eq: step 1}
    \matW_{\neq k}^\T \matW_{\neq k}^{}=\sum_{i=1}^{n_{-k}}\matW_{\neq k}(:,i)\matW_{\neq k}(:,i)^\T=\sum_{\veci_{-k}}\Tilde{\vecw}_k(\veci_{-k})\Tilde{\vecw}_k(\veci_{-k})^\T,
\end{equation}
where $\Tilde{\vecw}_k(\veci_{-k}):=\rmvec((\prod_{j=k+1}^{d}\matu_j(i_j)\prod_{j=1}^{k-1}\matu_j(i_j))^\T)$ and  $\veci_{-k}:=(i_{k+1}, \dots, i_d,$ $i_1, \dots, i_{k-1})\in[n_{k+1}]\times\cdots\times[n_d]\times[n_1]\times\cdots\times[n_{k-1}]$. 
By using $\rmvec(\matc\matx\matb^\T)=(\matb\otimes\matc)\rmvec(\matx)$ for matrices $\matB,\matC,\matx$ in appropriate size, we have
    \begin{align}
        \Tilde{\vecw}_k(\veci_{-k})&=
        \left(\left(\prod_{j=k+1}^{d}\matU_{j}(i_{j})\prod_{j=1}^{k-1}\matU_{j}(i_{j})\right)\otimes\matI_{r_k}\right)\rmvec(\matu_{k-1}(i_{k-1})^\T)\nonumber\\
        &=\prod_{j=k+1}^{d}(\matU_{j}(i_{j})\otimes\matI_{r_k})\prod_{j=1}^{k-1}(\matU_{j}(i_{j})\otimes\matI_{r_k})\rmvec(\matu_{k-1}(i_{k-1})^\T),\label{eq: step 2}
    \end{align}
where $\otimes$ denotes the Kronecker product.
Taking~\eqref{eq: step 2} into~\eqref{eq: step 1}, we can compute $\matW_{\neq k}^\T \matW_{\neq k}^{}$ in a recursive way, 
\begin{equation*}
    \begin{aligned}
        \tilde{\matH}_1&:=\sum_{i_{k-1}=1}^{n_{k-1}}\rmvec(\matu_{k-1}(i_{k-1})^\T)\rmvec(\matu_{k-1}(i_{k-1})^\T)^\T,\\
        \tilde{\matH}_2&:=\sum_{i_{k-2}=1}^{n_{k-2}}(\matU_{k-2}(i_{k-2})\otimes\matI_{r_k})\tilde{\matH}_1(\matU_{k-2}(i_{k-2})^\T\otimes\matI_{r_k}),\\
        &\ \vdots\\
        \tilde{\matH}_{d-1}&:=\sum_{i_{k+1}=1}^{n_{k+1}}(\matU_{k+1}(i_{k+1})\otimes\matI_{r_k})\tilde{\matH}_{d-2}(\matU_{k+1}(i_{k+1})^\T\otimes\matI_{r_k}).
    \end{aligned}
\end{equation*}
It follows that~$\matW_{\neq k}^\T \matW_{\neq k}^{}=\tilde{\matH}_{d-1}$. The computation of $\matW_{\neq k}^\T \matW_{\neq k}^{}$ for $k\in[d]$ requires
\[\sum_{k=1}^d\left(2n_{k-1}(r_kr_{k-1})^2 +2\sum_{i=1,i\neq k,k-1}^d n_ir_k^2r_ir_{i+1}(r_i+r_{i+1})\right)\] %The computation of $\matW_{\neq 1}^\T \matW_{\neq 1}^{}, \dots, \matW_{\neq d}^\T \matW_{\neq d}^{}$ requires
% $\sum_{k=1}^d\left(2n_{k-1}(r_kr_{k-1})^2 +2\sum_{i=1,i\neq k,k-1}^d n_ir_k^2r_ir_{i+1}(r_i+r_{i+1})\right)$
flops in Algorithm~\ref{alg: compute uneqk t uneqk}. If $n_1=\cdots=n_d$ and $r_1=\cdots=r_d$, computing the Riemannian gradient requires \[2d(d-1)\lvert \Omega \rvert r^3+2\lvert \Omega \rvert r^2+4d(d-2)nr^5+2dnr^4+C_{\mathrm{chol}}dr^6\] flops in total, which has the terms with order of $\mathcal{O}(n)$. In order to verify the improvement of Algorithm~\ref{alg: compute uneqk t uneqk}, we report a numerical comparison on a real dataset in Appendix~\ref{app: speedup}.

\begin{algorithm}[H]
    \caption{Efficient computation of $\matW_{\neq k}^\T \matW_{\neq k}^{}$}
    \label{alg: compute uneqk t uneqk}
    \begin{algorithmic}[1]
        \REQUIRE $k\in[d]$, core tensors $\tensU_l=\ten_{(2)}(\matW_l)$, slice matrices $\matu_l(i_l)$, $i_l\in[n_l],\ l\in[d]$.\STATE Set $j=\mathrm{mod}(k-2+d,d)+1$.
        \STATE Compute $\tilde{\matH}_1=\sum_{i_{j}=1}^{n_{j}}\rmvec(\matu_{j}(i_{j})^\T)\rmvec(\matu_{j}(i_{j})^\T)^\T$.
        \FOR{$l=2,3,\dots,d-1$}
            \STATE Set $j=\mathrm{mod}(k-l-1+d,d)+1$.
            \STATE Compute $\tilde{\matH}_l=\sum_{i_{j}=1}^{n_{j}}(\matU_{j}(i_{j})\otimes\matI_{r_k})\tilde{\matH}_{l-1}(\matU_{j}(i_{j})^\T\otimes\matI_{r_k})$.
        \ENDFOR
        \ENSURE $\tilde{\matH}_{d-1}=\matW_{\neq k}^\T \matW_{\neq k}^{}$.
    \end{algorithmic}
\end{algorithm}

\section{Convergence Analysis}\label{sec: convergence} 
The global convergence of TR-RGD and TR-RCG is analyzed in this section. Let $\{\vecmatW^{(t)}\}_{t\geq  0}$ be an infinite sequence generated by Algorithm~\ref{alg: TR-RGD} or Algorithm~\ref{alg: TR-RCG}. We prove that every accumulation point of $\{\vecmatW^{(t)}\}_{t\geq  0}$ is a stationary point; see Theorem~\ref{thm: global convergence} and Theorem~\ref{thm: global RCG}. 

\begin{lemma}[{\cite[Lemma 2.7]{boumal2019global}}]\label{lem: absil}
    Let $\tensM^\prime\subseteq\tensM$ be a compact Riemannian submanifold. Let $\mathcal{R}_x:\tangent_x\tensM^\prime\to\tensM^\prime$ be retraction. $f:\tensM^\prime\to\mathbb{R}$ has Lipschitz continuous gradient on the convex hull of  $\tensM^\prime$. Then, there exists a constant $L>0$, such that 
    \[\left\lvert f(\mathcal{R}_x(\xi))-f(x)-g_x(\xi,\grad\!f(x))\right\rvert\leq \frac{L}{2} g_x(\xi,\xi)\quad\text{for all}\ x\in\tensM^\prime,\ \xi\in \tangent_x\tensM^\prime.\]
\end{lemma}

Since the search space $\tensM$ is flat, the retraction map in Lemma~\ref{lem: absil} is chosen as the identity map. We observe the coercivity of $f$ in~\eqref{eq: newly stated problem} from the regularization term~$\frac{\lambda}{2}\|\vecmatW\|_\mathrm{F}^2$. Hence, the level set $\tensL:=\{\vecmatW:\ f(\vecmatW)\leq  f(\vecmatW^{(0)})\}$ is compact. The sequence of function values $\{f(\vecmatW^{(t)})\}_{t\geq  0}$, obtained from~\eqref{eq: linmin} and~\eqref{eq: Armijo}, is monotonically decreasing. It holds that $\{\vecmatW^{(t)}\}_{t\geq  0}$ is a bounded sequence in $\tensL$ with
{\[\|\vecmatW^{(t)}\|^2_\mathrm{F}=\frac{2\left(f(\vecmatW^{(t)})-f_\Omega(\vecmatW^{(t)})\right)}{\lambda}\leq\frac{2f(\vecmatW^{(t)})}{\lambda} \leq \frac{2f(\vecmatW^{(0)})}{\lambda}.\]}Therefore, there exist accumulation points for the sequence $\{\vecmatW^{(t)}\}_{t\geq  0}$. Moreover, the objective function $f$ has Lipschitz continuous gradient. By using Lemma~\ref{lem: absil} and~\cite[Proposition~4.3]{dong2022new}, we can prove the global convergence of the proposed Riemannian gradient descent algorithm in a same fashion. 

\begin{theorem}\label{thm: global convergence}
    Let $\{\vecmatW^{(t)}\}_{t\geq  0}$ be an infinite sequence generated by Algorithm~\ref{alg: TR-RGD}. Then, there exists $C>0$ such that $f(\vecmatW^{(t)})-f(\vecmatW^{(t+1)})>C\|\grad\!f(\vecmatW^{(t)})\|_{\vecmatW^{(t)}}$. Furthermore, we have: 1) every accumulation point of $\{\vecmatW^{(t)}\}_{t\geq  0}$ is a stationary point of $f$; 2) the algorithm returns $\vecmatW\in\tensM$ satisfying $\|\grad\!f(\vecmatW)\|_{\vecmatW}<\epsilon$ after $\left\lceil{f(\vecmatW^{(0)})}/{(C\epsilon^2)}\right\rceil$ iterations at most.
\end{theorem}

Now, we discuss the global convergence of the proposed Riemannian conjugate gradient algorithm in Algorithm~\ref{alg: TR-RCG}; interested readers are referred to~\cite{absil2009optimization,sato2022riemannian} for the convergence of RCG on general manifolds. Here, we follow the convergence analysis in~\cite{absil2009optimization}. To fulfill the basic assumptions in ~\cite[Theorem 4.3.1]{absil2009optimization}, i.e., the sequence of search directions $\{\vecmateta^{(t)}\}_{t\geq  0}$ is \emph{gradient-related} to $\{\vecmatW^{(t)}\}_{t\geq  0}$, we enforce $\vecmateta^{(t)}$ to be 
\begin{equation}
    {\vecmateta^{(t)}=-\grad\!f(\vecmatW^{(t)})\quad \text{if}\ g_{\vecmatW^{(t)}}(\vecmateta^{(t-1)},\grad\!f(\vecmatW^{(t)}))\geq 0}.\label{eq: restart}
\end{equation}
Therefore, $\{\vecmateta^{(t)}\}_{t\geq  0}$ are descent directions with
\begin{align}
    &g_{\vecmatW^{(t)}}(\vecmateta^{(t)},\grad\!f(\vecmatW^{(t)}))\nonumber\\
    =&-\|\grad\!f(\vecmatW^{(t)})\|_{\vecmatW^{(t)}}^2\nonumber+\beta^{(t)}g_{\vecmatW^{(t)}}(\vecmateta^{(t-1)},\grad\!f(\vecmatW^{(t)}))\nonumber\\
    \leq&-\|\grad\!f(\vecmatW^{(t)})\|_{\vecmatW^{(t)}}^2.\label{eq: fact 1}
\end{align}
In addition, it follows from Algorithm~\ref{alg: TR-RCG} that 
\begin{equation}
    \|\vecmateta^{(t)}\|_\mathrm{F}=\frac{\|\vecmatW^{(t+1)}-\vecmatW^{(t)}\|_\mathrm{F}}{s^{(t)}}\leq\frac{\|\vecmatW^{(t+1)}\|_\mathrm{F}+\|\vecmatW^{(t)}\|_\mathrm{F}}{s_{\min}}\leq\frac{2}{s_{\min}}\sqrt{\frac{2f(\vecmatW^{(0)})}{\lambda}}\label{eq: fact 2}
\end{equation} is uniformly bounded. In view of \eqref{eq: fact 1} and~\eqref{eq: fact 2}, $\{\vecmateta^{(t)}\}_{t\geq  0}$ is gradient-related to $\{\vecmatW^{(t)}\}_{t\geq  0}$. By using~\cite[Theorem 4.3.1]{absil2009optimization}, we have the following result. 
\begin{theorem}\label{thm: global RCG}
    Let $\{\vecmatW^{(t)}\}_{t\geq 0}$ be an infinite sequence generated by Algorithm~\ref{alg: TR-RCG} with modified search direction~\eqref{eq: restart}. Then, every accumulation point of $\{\vecmatW^{(t)}\}_{t\geq 0}$ is a stationary point of $f$. 
\end{theorem}

\section{Numerical Experiments}\label{sec: numerical exps}
In this section, we numerically compare TR-RGD (Algorithm~\ref{alg: TR-RGD}) and TR-RCG (Algorithm~\ref{alg: TR-RCG}) with existing algorithms based on different tensor decompositions on synthetic and real-world datasets, including movie ratings, hyperspectral images, and high-dimensional functions. First, we introduce all the compared algorithms and default settings.

We consider the Riemannian conjugate gradient algorithm\footnote{TTeMPS toolbox: \url{https://www.epfl.ch/labs/anchp/index-html/software/ttemps/}.} in~\cite{steinlechner2016riemannian} based on tensor train decomposition, which is denoted by ``TT-RCG". 
For CP-based algorithm, we choose ``{CP-WOPT}"~\cite{acar2011scalable} in {Tensor-Toolbox\footnote{Tensor-Toolbox v3.4: \url{http://www.tensortoolbox.org/}.}} 
by Bader and Kolda~\cite{bader2008efficient} with limited-memory BFGS algorithm\footnote{Available from \url{https://github.com/stephenbeckr/L-BFGS-B-C}.} recommended by the authors. 
``GeomCG"~\cite{kressner2014low} is a Riemannian conjugate gradient algorithm\footnote{GeomCG toolbox: \url{https://www.epfl.ch/labs/anchp/index-html/software/geomcg/}.} for Tucker-based tensor completion problem. 
In addition, we consider a nuclear-norm-based algorithm\footnote{Available from \url{https://github.com/andrewssobral/mctc4bmi/tree/master/algs_tc/LRTC}.} ``HaLRTC"~\cite{liu2012tensor,sobral_sbmi_prl_2016}. If not specified, we adopt default settings for all the compared algorithm. 

Note that there are different ways to choose stepsize. In the preliminary numerical experiments, TR-RGD with stepsize RBB2  in~\eqref{eq: alternating BB stepsizes} performs better than the algorithm with RBB1. Therefore, we only consider stepsize RBB2  in the following comparisons. ``TR-RGD (exact)" and ``TR-RGD (RBB)" denote Algorithm~\ref{alg: TR-RGD} with exact line search~\eqref{eq: linmin} and Algorithm~\ref{alg: TR-RGD} with Armijo backtracking~\eqref{eq: Armijo} and RBB2, respectively. Moreover, the vanilla Euclidean gradient descent algorithm, denoted by ``TR-GD", is implemented for comparison. The stepsize for TR-GD is based on Armijo backtracking and the standard BB~\cite{bb1988}. The default settings of line search parameters are {$\rho=0.4$, $a=10^{-5}$, and $s_{\min}=10^{-10}$}. In addition, we implement the alternating least squares algorithm~\cite{wang2017efficient,zhao2019learning} called ``TR-ALS" for the tensor completion problem~\eqref{eq: newly stated problem}.

All algorithms are initialized from $\tensX^{(0)}=\tau(\vecmatW)$ in which $\tau$ is defined in~\eqref{eq: mapping tau} and $\vecmatW$ is randomly generated. {We define the relative error on a sampling set $\Omega$
\[\varepsilon_{\Omega}(\vecmatW):=\frac{\|\proj_\Omega(\tau(\vecmatW))-\proj_\Omega(\tensA)\|_\mathrm{F}}{\|\proj_\Omega(\tensA)\|_\mathrm{F}}.\] 
We refer to the training error as $\varepsilon_\Omega(\vecmatW)$. Furthermore, we evaluate the test error $\varepsilon_\Gamma(\vecmatW)$ on a test set $\Gamma$ different from $\Omega$. The default setting of $\lvert \Gamma \rvert $ is $100$.} For stopping criteria, we terminate the algorithms once one of the following criteria is reached: 1) the relative error $\varepsilon_{\Omega}(\vecmatW^{(t)})<10^{-12}$; 2) the relative change $\lvert{(\varepsilon_{\Omega}(\vecmatW^{(t)})-\varepsilon_{\Omega}(\vecmatW^{(t-1)}))}/{\varepsilon_{\Omega}(\vecmatW^{(t-1)})}\rvert<\varepsilon$; 3)~the gradient norm $\|\grad\!f(\vecmatW)\|_{\mathrm{F}}<\varepsilon$. This criterion is only activated for Riemannian methods; 4) maximum iteration number; 5) time budget. The tolerance $\varepsilon$ is chosen as $10^{-8}$. All experiments are performed on a MacBook Pro 2019 with MacOS Ventura 13.1, 2.4 GHz 8 core Intel Core i9 processor, 32GB memory, and Matlab R2020b. The codes of TR-RGD and TR-RCG can be downloaded from~\url{https://github.com/JimmyPeng1998}.

\subsection{Synthetic data}\label{subsec: synthetic}
In this subsection, we investigate the numerical performance and the reconstruction ability of tensor completion algorithms in TR decomposition on noiseless and noisy observations. 

\paragraph{Noiseless observations}
We consider a synthetic low-rank tensor $\tensA\in\mathbb{R}^{n_1\times \cdots \times n_d}$ in~\eqref{eq: newly stated problem} generated by \[\tensA=\tau(\vecmatW),\] where $\tau$ is defined in~\eqref{eq: mapping tau}, each entry of $\vecmatW\in\tensM$ is uniformly sampled from $[0,1]$, $d=3$, $n_1=n_2=n_3=100$, and TR rank $\vecr^*=(6,6,6)$. Given the sampling rate $p$, we formulate the sampling set $\Omega$ by randomly selecting $pn_1\cdots n_d$ samples from $[n_1]\times\cdots\times[n_d]$. In order to get an unbiased recovery result, we choose the regularization parameter $\lambda=0$, the sampling rate $p=0.3$, and $\vecr=(6,6,6)$. The time budget is $1800$s and the maximum iteration number is $10000$.

\begin{figure}[htbp]
	\centering
	\subfigure{\includegraphics[width=0.45\textwidth]{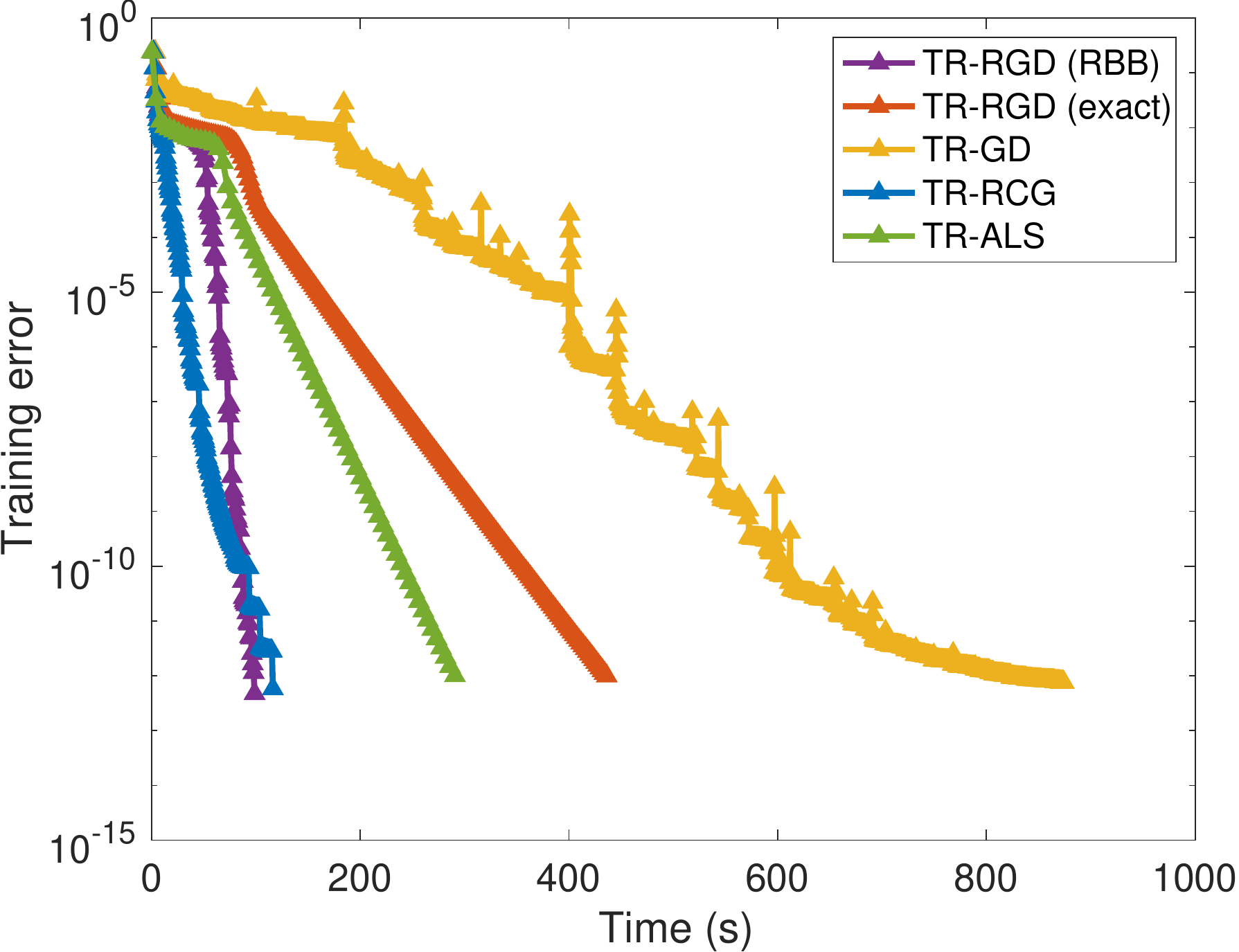}}\quad
	\subfigure{\includegraphics[width=0.45\textwidth]{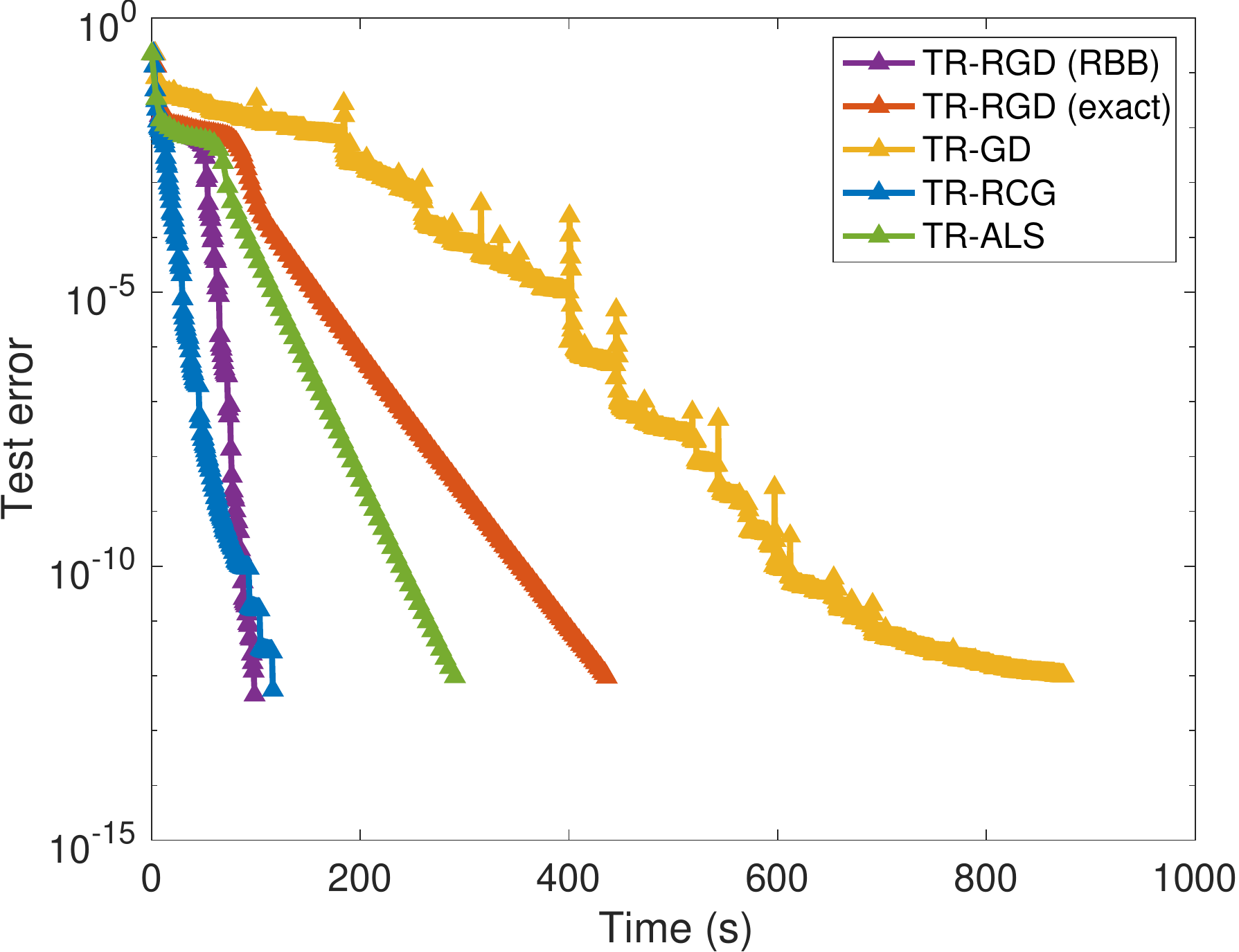}}
	\caption{Numerical results of noiseless case. Left: training error. Right: test error}
	\label{fig: noiseless final results}
\end{figure}

Figure~\ref{fig: noiseless final results} presents the numerical results of noiseless case. First, we observe that all algorithms successfully recover the true tensor within time budget. TR-RGD (RBB) and TR-RCG perform better than the alternating least squares algorithm (TR-ALS) in terms of training error and test error, and TR-RGD (RBB) is comparable to TR-RCG. Second, TR-RGD with RBB stepsize is more efficient than the algorithm with exact line search, since computing exact line seach needs to find the roots for a polynomial of degree $2d-1$. However, computing the coefficients of such polynomial is expensive in practice. Third, it is worth noting that the proposed Riemannian gradient algorithms outperform the Euclidean gradient algorithm with BB stepsizes, implying that the new metric have a preconditioning effect indeed.

\begin{figure}[htbp]
    \centering
    \subfigure{\includegraphics[width=0.32\textwidth]{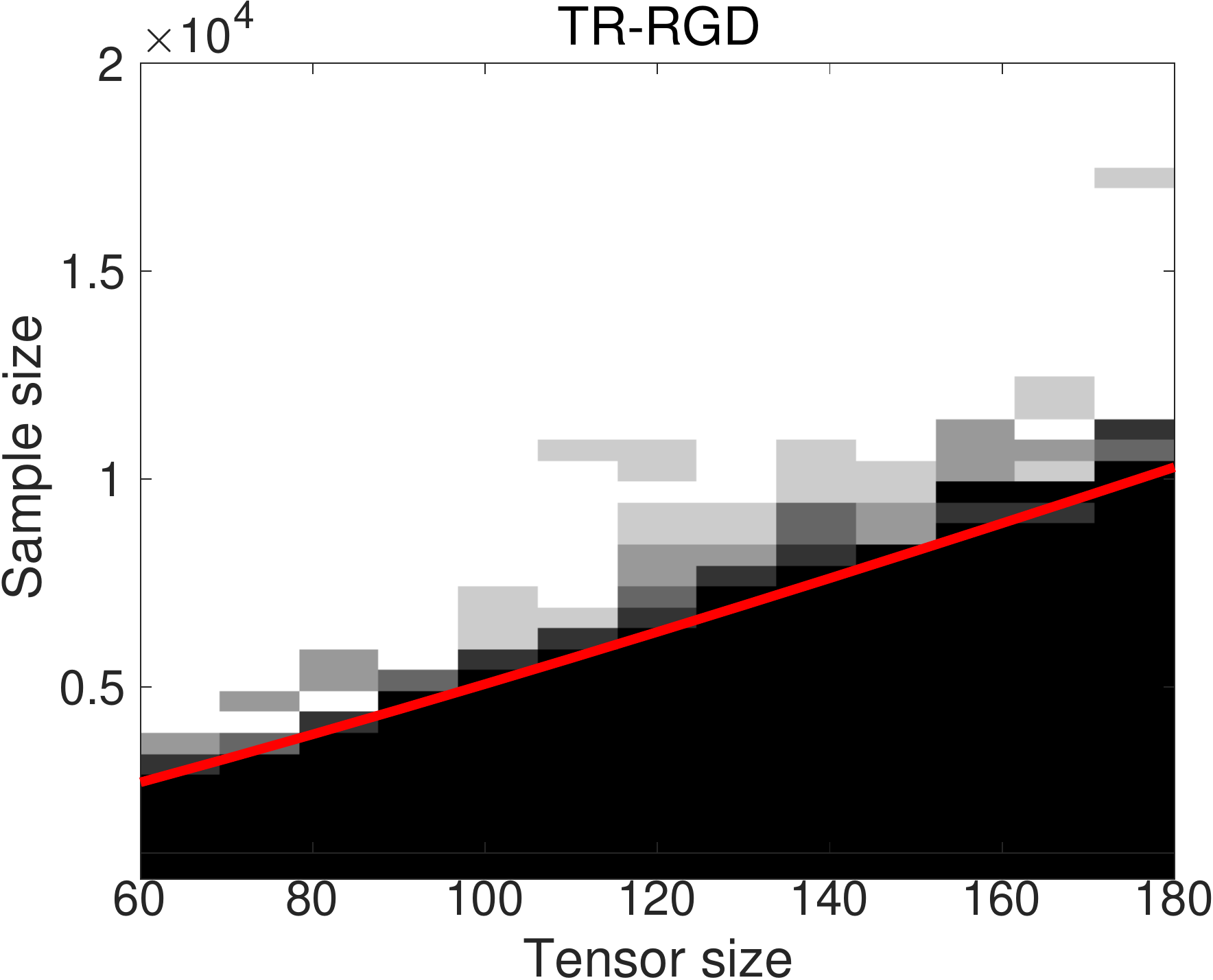}}
    \subfigure{\includegraphics[width=0.32\textwidth]{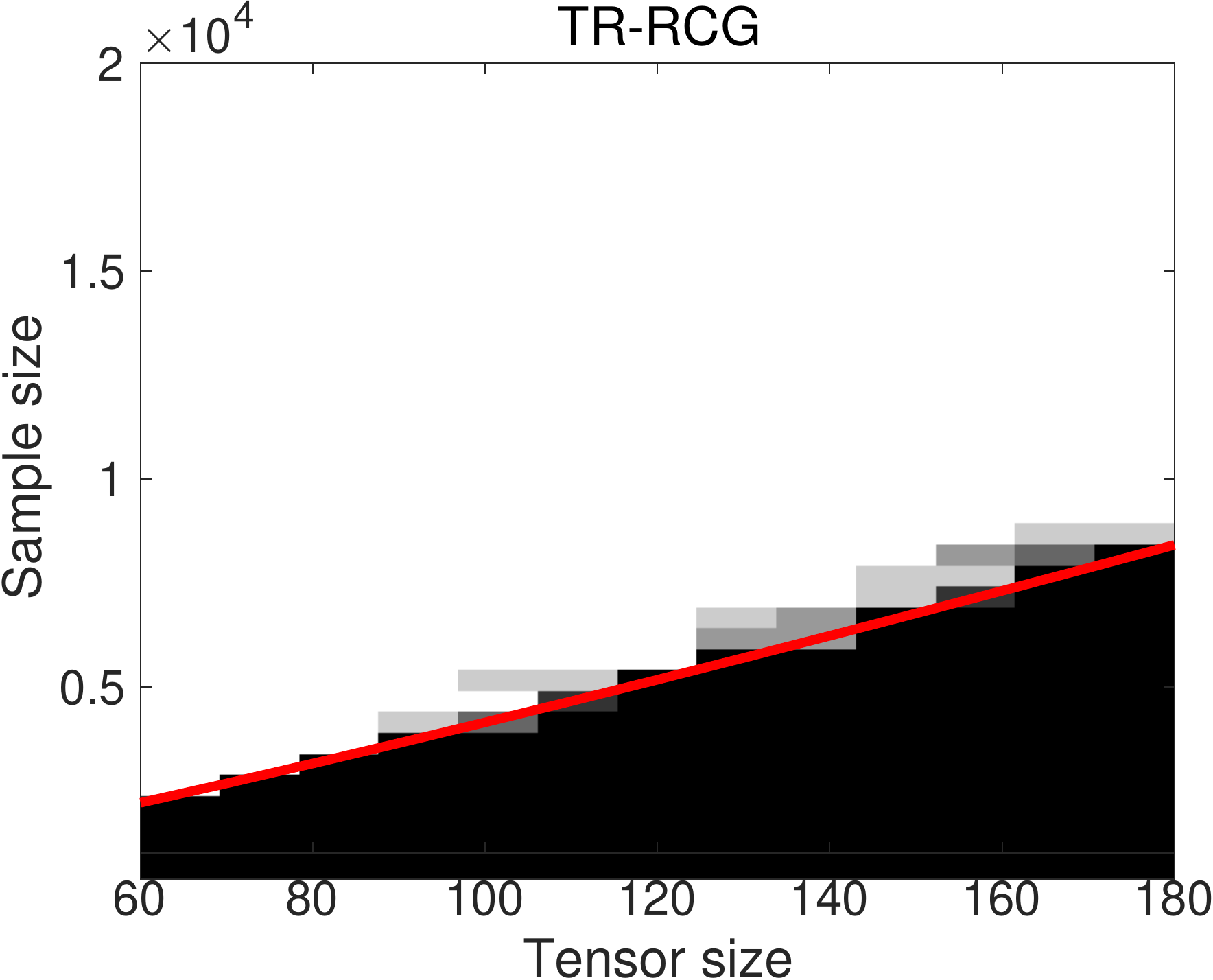}}
    \subfigure{\includegraphics[width=0.32\textwidth]{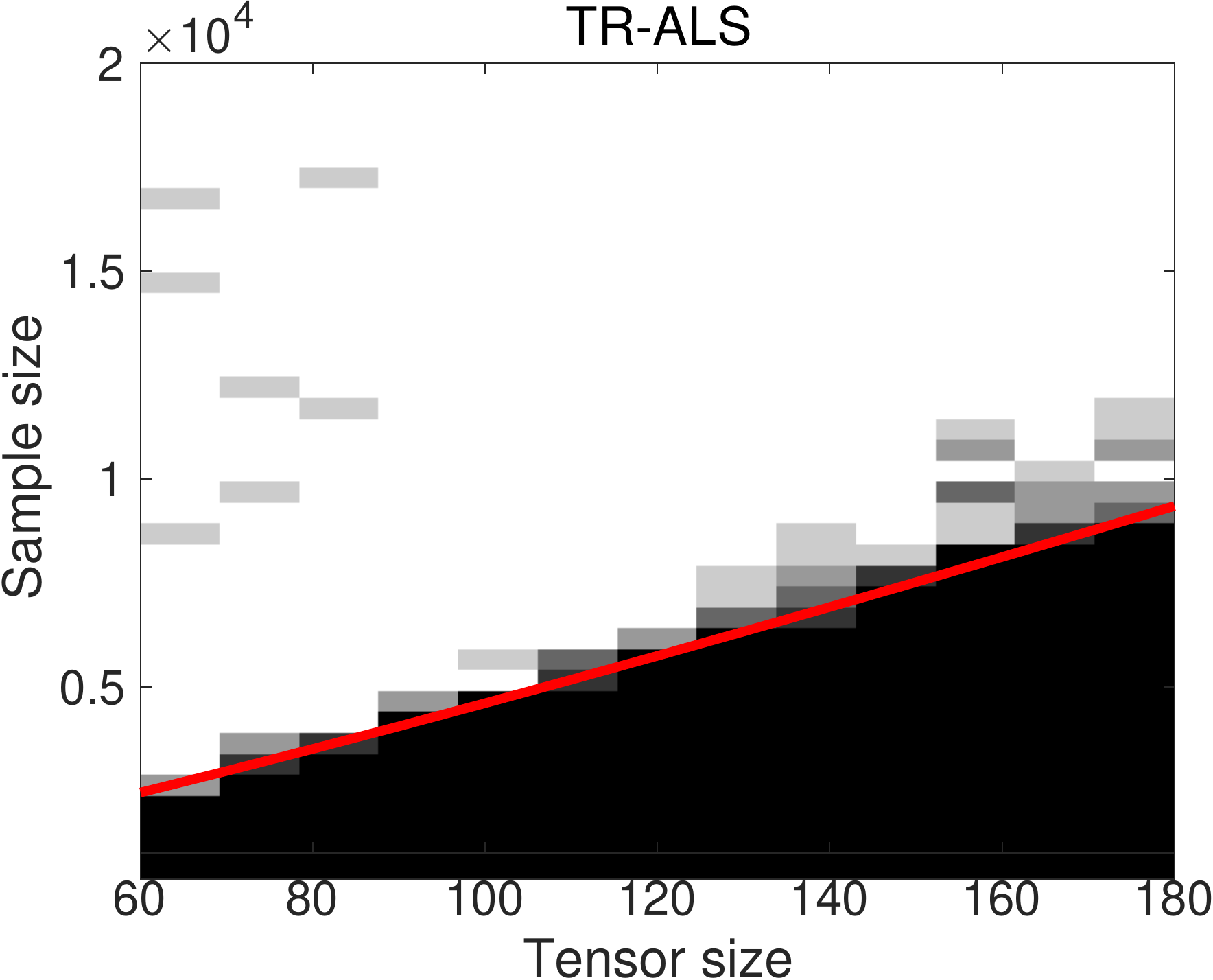}}
    \caption{Phase plots of recovery results for five runs. The white block indicates successful recovery in all five runs, while the black block signifies failure of recovery in all five runs}
    \label{fig: phase plots}
\end{figure}

The question ``how many samples are required to recover a low-rank data'' is interesting but challenging. In the matrix case $(d=2)$, around $\mathcal{O}(nr\log n)$ samples are necessary to recover a low-rank matrix~\cite{keshavan2009matrix,candes2010power}. We aim to numerically investigate this question by exploring the relationship between sampling rate and tensor size $n$ for third-order tensors with TR rank $\vecr^*=(3,3,3)$. To this end, we randomly generate third order synthetic tensors $\tensA$ with TR rank $\mathbf{r}=(3,3,3)$ in the same fashion as above, and select tensor size $n:=n_1=n_2=n_3$ from $\{60,70,\dots,180\}$ and the sample size $\lvert \Omega \rvert$ from $\{1000,1500,\dots,20000\}$. We run TR-RGD, TR-RCG, and TR-ALS algorithms five times for each combination of $n$ and $\lvert \Omega \rvert$. One algorithm is regarded to successfully recover a tensor if the test error $\varepsilon_{\Gamma}<10^{-4}$ within maximum iteration number 250. Figure~\ref{fig: phase plots} illustrates the phase plots of recovery results, in which the gray level of each block represents the number of successful recovery. White blocks imply successful recovery in all five runs. The red line represents $\mathcal{O}(n\log n)$ when $n$ scales. The phase plots suggest a similar behavior to the matrix case, which is consistent to other numerical experiments for $d=3$, e.g.,~\cite{kressner2014low}.

\paragraph{Noisy observations}
Furthermore, we investigate the reconstruction ability of TR-based algorithms under different noise levels. We consider a synthetic noisy tensor 
\begin{equation*}
    \tensA{:=}\frac{\hat{\tensA}}{\|\hat{\tensA}\|_\mathrm{F}}+\sigma\cdot\frac{\tensE}{\|\tensE\|_\mathrm{F}}, 
\end{equation*}
where $\hat{\tensA}$ is generated in the same rule as the noiseless case with TR rank $\mathbf{r}^*=(3,3,3)$ and $n_1=n_2=n_3=100$. $\tensE$ is a tensor with its entries being sampled from normal distribution $\tensN(0,1)$. The parameter $\sigma$ measures the noise level of a tensor. Ideally, the relative errors $\varepsilon_\Omega$ and $\varepsilon_\Gamma$ are supposed to be the noise level~$\sigma$ when the algorithm terminates. We set $\sigma=10^{-3},10^{-4},10^{-5},10^{-6},10^{-7},10^{-8}$, and the regularization parameter $\lambda=10^{-12}$ to get unbiased recovery results. We observe from Fig.~\ref{fig: phase plots} that it requires at least $8000$ samples to recover the noiseless tensor, i.e., $p$ should be larger than $8000/100^3=0.008$. Therefore, we choose the sampling rate $p=0.05$. The time budget is $120$s and the maximum iteration number is $1000$. 

\begin{figure}[htbp]
    \centering
    \subfigure{\includegraphics[height=3.65cm]{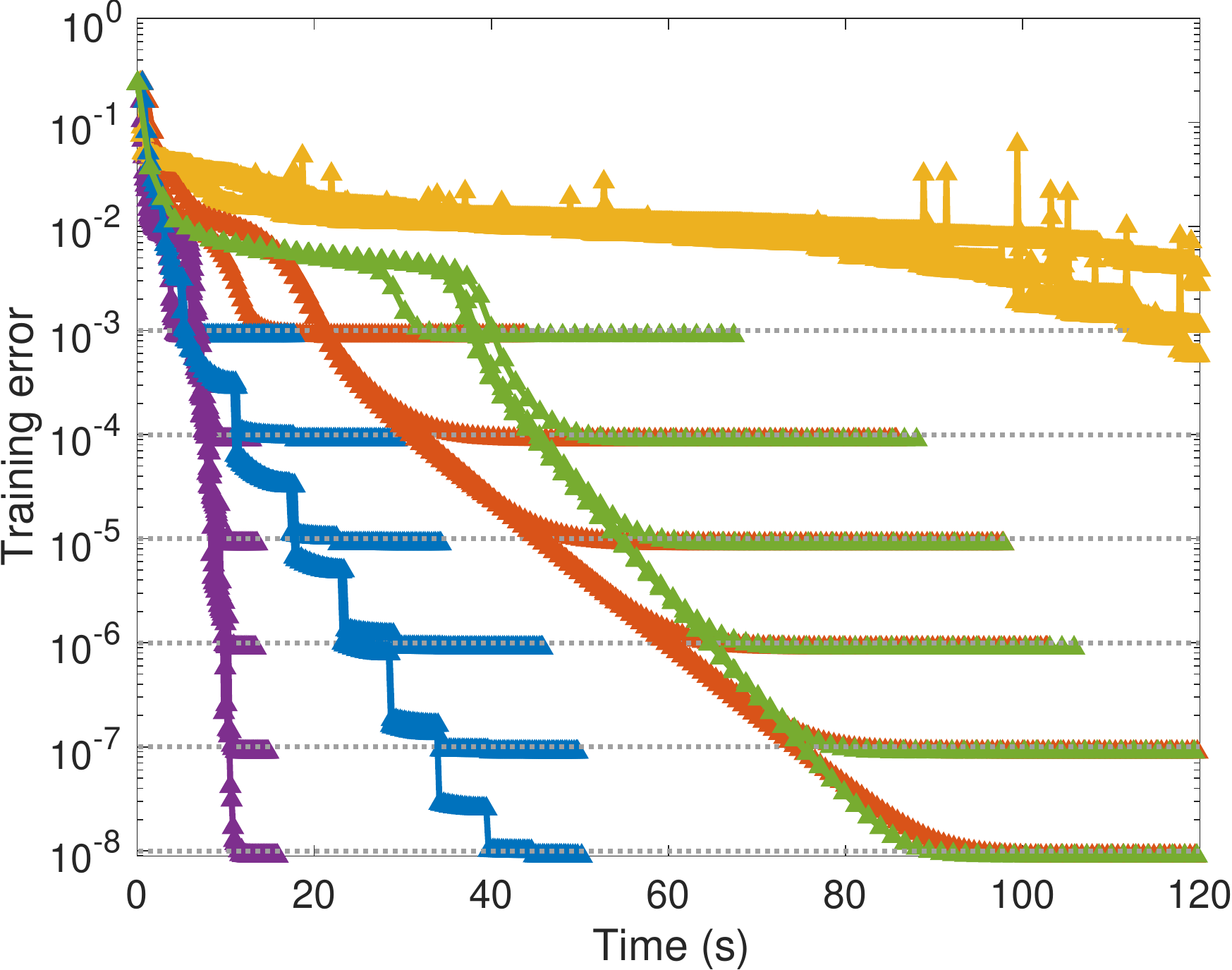}}\qquad
    \subfigure{\includegraphics[height=3.9cm]{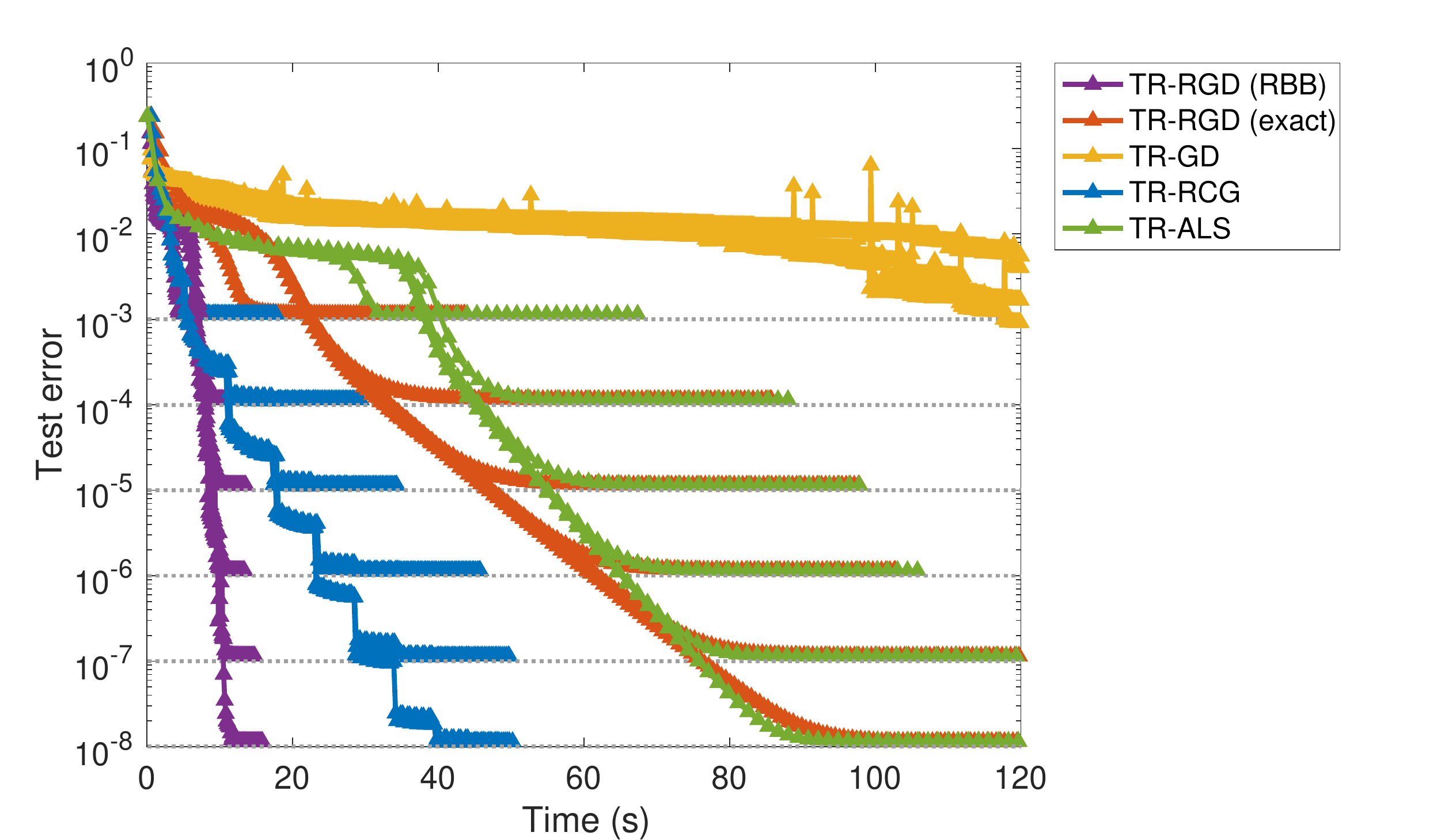}}
    \caption{Reconstruction ability of TR-based algorithms under different noise levels. Left: training error. Right: test error}
    \label{fig: Noisy reconstruction}
\end{figure}

Figure~\ref{fig: Noisy reconstruction} shows the recovery performance of TR-based algorithms under different noise levels. All algorithms successfully recover the underlying low-rank tensor within time budget except TR-GD. Since the TR rank parameter $\vecr$ is exactly the true TR rank $\vecr^*$ of data tensor $\tensA$, we observe from Table~\ref{tab: results for noisy} that the test errors are comparable to and slightly larger than the training errors.

\begin{table}[htbp]
    \centering
    \caption{Training and test errors of noisy case}
    \scriptsize
    \resizebox{\textwidth}{!}{
    \begin{tabular}{ccccccc}
        \toprule
        $\sigma$ & Error & TR-RGD (RBB) & TR-RGD (exact) & TR-RCG & TR-GD & TR-ALS \\
        \midrule
        \multirow{2}*{1e-3} & Training & 8.9142e-04 & 8.9142e-04 & 8.9142e-04 & 1.1003e-03 & 8.9142e-04 \\
         & Test & 1.1471e-03 & 1.1471e-03 & 1.1472e-03 & 1.4066e-03 & 1.1470e-03 \\
        % \midrule
        \multirow{2}*{1e-4} & Training & 8.9154e-05 & 8.9154e-05 & 8.9154e-05 & 1.1129e-03 & 8.9154e-05 \\ 
         & Test & 1.1458e-04 & 1.1459e-04 & 1.1461e-04 & 1.6735e-03 & 1.1458e-04 \\ 
        % \midrule
        \multirow{2}*{1e-5} & Training & 8.9155e-06 & 8.9155e-06 & 8.9155e-06 & 5.7432e-04 & 8.9155e-06 \\ 
         & Test & 1.1457e-05 & 1.1457e-05 & 1.1461e-05 & 9.1646e-04 & 1.1457e-05 \\ 
        % \midrule
        \multirow{2}*{1e-6} & Training & 8.9155e-07 & 8.9155e-07 & 8.9155e-07 & 2.7753e-03 & 8.9155e-07 \\ 
         & Test & 1.1457e-06 & 1.1457e-06 & 1.1458e-06 & 4.0237e-03 & 1.1458e-06 \\ 
        % \midrule
        \multirow{2}*{1e-7} & Training & 8.9155e-08 & 8.9155e-08 & 8.9155e-08 & 3.8170e-03 & 8.9156e-08 \\ 
         & Test & 1.1456e-07 & 1.1456e-07 & 1.1457e-07 & 5.4727e-03 & 1.1462e-07 \\ 
        % \midrule
        \multirow{2}*{1e-8} & Training & 8.9156e-09 & 8.9156e-09 & 8.9156e-09 & 9.5924e-04 & 8.9289e-09 \\ 
         & Test & 1.1450e-08 & 1.1455e-08 & 1.1451e-08 & 1.4597e-03 & 1.1512e-08 \\ 
        \bottomrule
    \end{tabular}
    }
    \label{tab: results for noisy}
\end{table}

\subsection{Experiments on MovieLens 1M dataset}
We consider a real-world tensor completion problem on the {MovieLens 1M\footnote{Available from \url{https://grouplens.org/datasets/movielens/1m/}.}} dataset including one million movie ratings from 6040 users on 3952 movies from September 19th, 1997 to April 22nd, 1998. We formulate the movie ratings as a third order tensor $\tensA$ of size $6040\times 3952\times 150$ by choosing one week as a period. We randomly select $80\%$ of the known ratings as training set $\Omega$ and the rest $20\%$ ratings are test set $\Gamma$. We compared the proposed algorithms TR-RGD and TR-RCG with other tensor completion algorithms including TR-GD, TR-ALS, TT-RCG, CP-WOPT, and geomCG. We choose the parameters: TR rank $\vecr=(6,10,3)$ and $\vecr=(6,6,6)$, regularization parameter $\lambda=1$, and the time budget 500s. To ensure a close number of parameters in different search spaces, we choose the TT rank as $(1,10,10,1)$ and $(1,9,9,1)$, the Tucker rank as $(60,30,18)$ and $(36,36,36)$, and the CP rank as $49$ and $36$; see Appendix~\ref{app: search space} for details.

\begin{figure}[htbp]
    \centering
    {\includegraphics[width=0.45\textwidth]{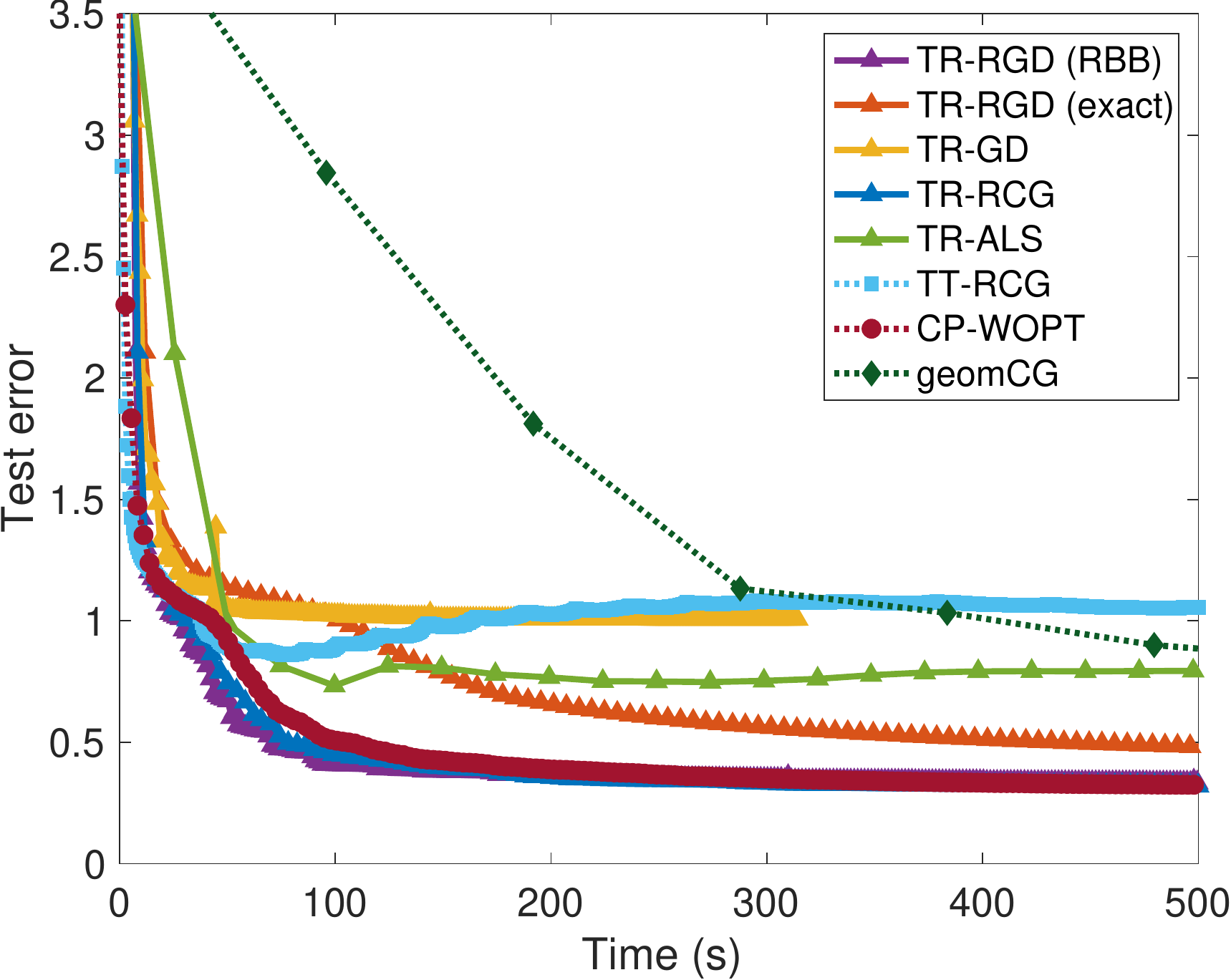}}\quad
    {\includegraphics[width=0.45\textwidth]{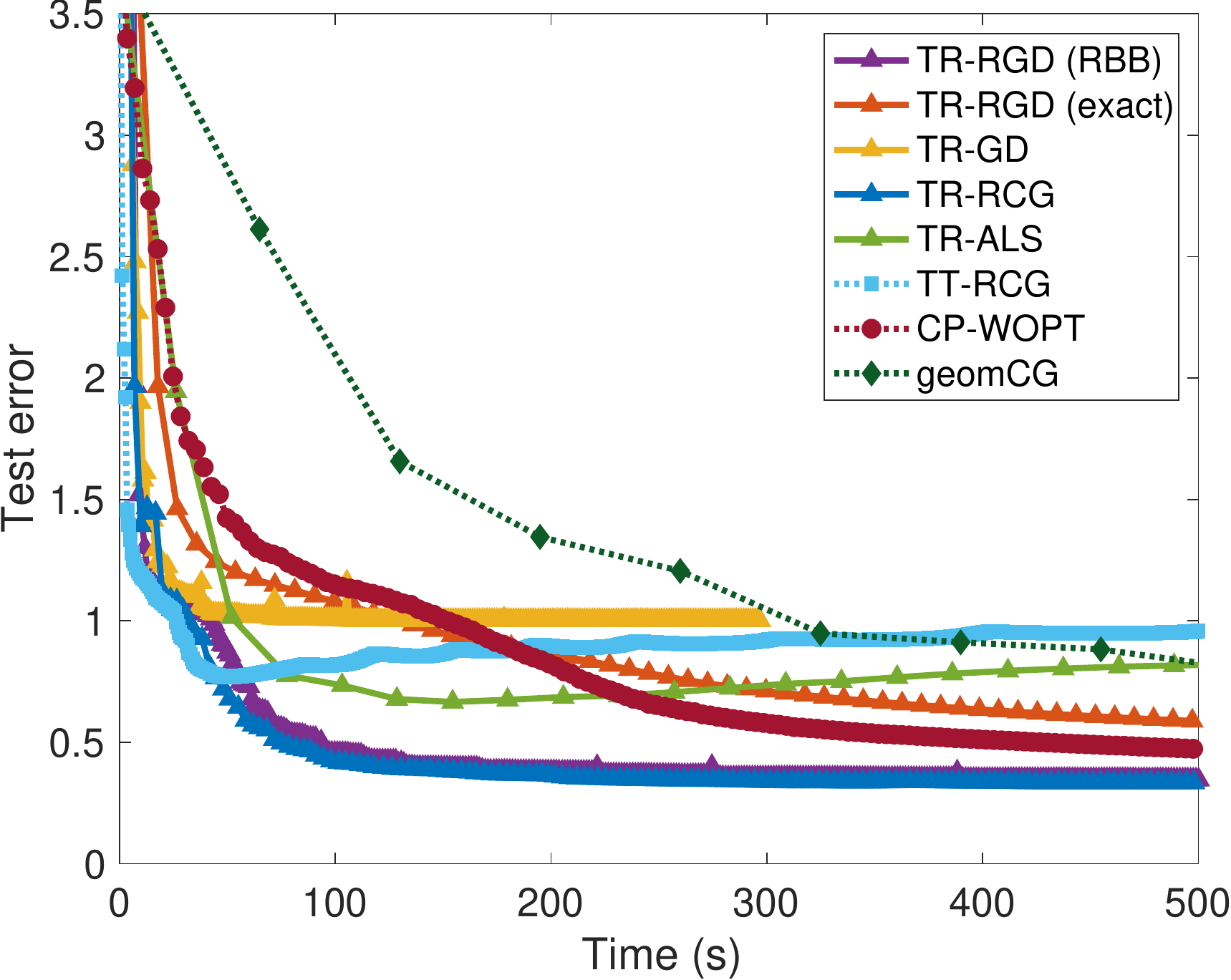}}
    \caption{Test error on MovieLens 1M dataset. Left: $\mathbf{r}=(6,6,6)$. Right: $\mathbf{r}=(6,10,3)$}
    \label{fig: test error}
\end{figure}

The test error under two different rank selections are shown in Fig.~\ref{fig: test error}. We observe that the proposed TR-RGD (RBB) and TR-RCG are comparable with CP-WOPT since CP-WOPT also benefits from the second-order information. They have faster convergence with lower test errors than the others. Moreover, TR-RGD (exact) spends more time than TR-RGD (RBB). For the sake of brevity, we only consider TR-RGD (RBB) in the following experiments and simplify ``TR-RGD (RBB)'' to ``TR-RGD''.

\subsection{Experiments on hyperspectral images}
\begin{figure}[htbp]
    \centering
    \subfigure{\includegraphics[height=4.2cm, trim=0 0.8 0 0]{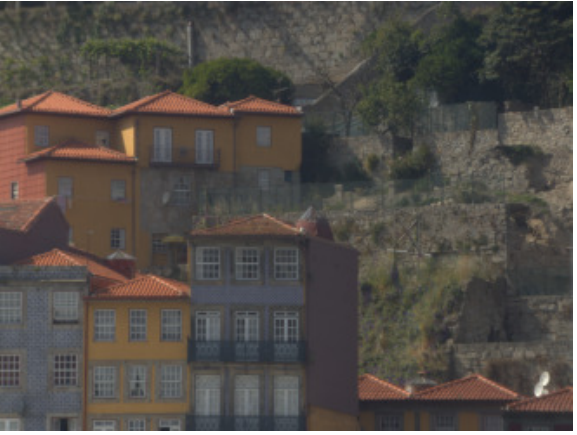}}
    \quad
    \subfigure{\includegraphics[height=4.16cm, trim=0 0 0 0.8]{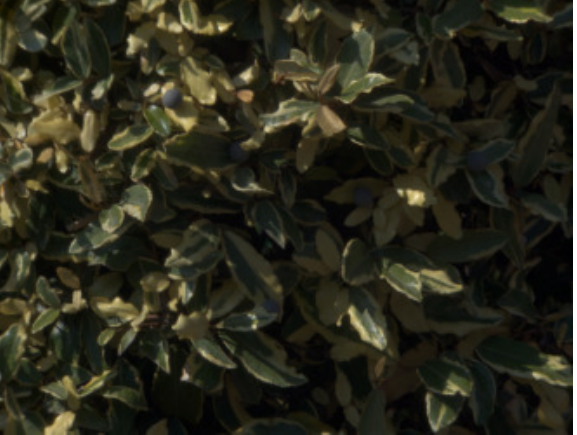}}
    \caption{Hyperspectral images. Left: ``Ribeira House Shrubs''. Right: ``Bom Jesus Bush''}
    \label{fig: Colored HSI}
\end{figure}

In this experiment, we consider the completion task on hyperspectral images. A hyperspectral image is formulated as a third order tensor $\tensA\in\mathbb{R}^{n_1\times n_2\times n_3}$. Mode three of $\tensA$ represents $n_3$ wavelength values of light. Mode one and two represents reflectance level of light under different wavelengths. We select two images\footnote{Image source: hsi\_34.mat and hsi\_46.mat from \url{https://figshare.manchester.ac.uk/articles/dataset/Fifty_hyperspectral_reflectance_images_of_outdoor_scenes/14877285}.} from {``50 reduced hyperspectral reflectance images''} by Foster~\cite{foster2022colour}: ``Ribeira Houses Shrubs'', abbreviated as ``Ribeira'', with size $249\times 330\times 33$; and ``Bom Jesus Bush'', abbreviated as ``Bush'', with size $250\times 330\times 33$. We transform hyperspectral images to RGB images by following the {tutorial\footnote{\url{https://personalpages.manchester.ac.uk/staff/david.foster/Tutorial_HSI2RGB/Tutorial_HSI2RGB.html}.}}. Figure~\ref{fig: Colored HSI} shows these two hyperspectral images represented as RGB images. 

To evaluate the recovery performance of image completion, 
peak signal-to-noise ratio (PSNR) is used to measure the similarity of two images, defined by
\[\mathrm{PSNR}:=10\log_{10}\left(\frac{\max(\tensA)^2}{\mathrm{MSE}}\right)=10\log_{10}\left(n_1n_2n_3\frac{\max(\tensA)^2}{\|\tensX-\tensA\|_\mathrm{F}^2}\right),\]
where $\max(\tensA)$ denotes the highest pixel value of $\tensA$, and $\mathrm{MSE}$ is the mean square error defined by $\mathrm{MSE}:={\|\tensX-\tensA\|_\mathrm{F}^2}/(n_1 n_2 n_3)$. The relative error \[\mathrm{relerr}(\tensX):=\frac{\|\tensX-\tensA\|_\mathrm{F}}{\|\tensA\|_\mathrm{F}}\] is also reported. In our experiments, the PSNR between a recovered image and original image will be reported after being transformed into RGB images. Given sampling rate $p$, we formulate the sampling set $\Omega$ in the same fashion as subsection~\ref{subsec: computational skills}. The initial guess $\tensX^{(0)}$ is randomly generated with normal distribution~$\mathcal{N}(0,1)$. We set the TR rank as $\vecr=(7,16,7)$, the Tucker rank as (65,65,7) following~\cite{kressner2014low}, the TT rank as $(1,15,15,1)$, and the CP rank as~$110$. The search spaces are of a similar size for different tensor formats; see Appendix~\ref{app: search space}. Moreover, the nuclear-norm-based algorithm, HaLRTC~\cite{liu2012tensor}, is compared. The maximum iteration number is 200.

\begin{figure}[htbp]
    \centering

    \begin{tikzpicture}
        \coordinate (O) at (0,0);
        \coordinate (C) at (1.7,0); % 控制图片中心横向间距
        \coordinate (R) at (0,-1.48); % 控制图片中心纵向间距
        \coordinate (d) at (0,0.52); % 控制db离图片中心的距离
        \coordinate (t) at (0,0.6); % 控制算法离图片中心的距离
        \node[inner sep=0] (image) at (0,0) {\includegraphics[width=1.69cm, height=1.25cm]{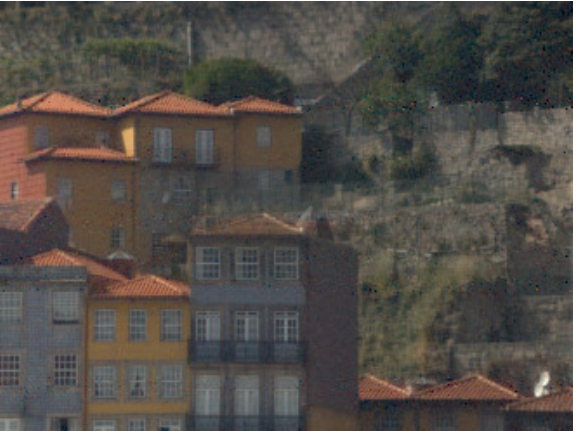}};
        \node[above] at ($(O)+(t)$) {\small TR-RCG};
        \node[below] at ($(O)-(d)$) {\scriptsize 38.4310dB};
        
        \node[inner sep=0] (img) at ($(O)+(C)+(0,0.005)$) {\includegraphics[width=1.69cm, height=1.25cm]{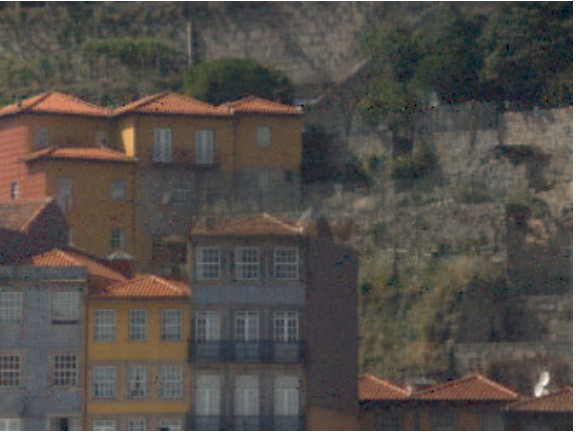}};
        \node[above] at ($(img)+(t)$) {\small TR-RGD};
        \node[below] at ($(img)-(d)$) {\scriptsize 37.5599dB};
        
        \node[inner sep=0] (img) at ($(O)+2*(C)$) {\includegraphics[width=1.69cm, height=1.25cm]{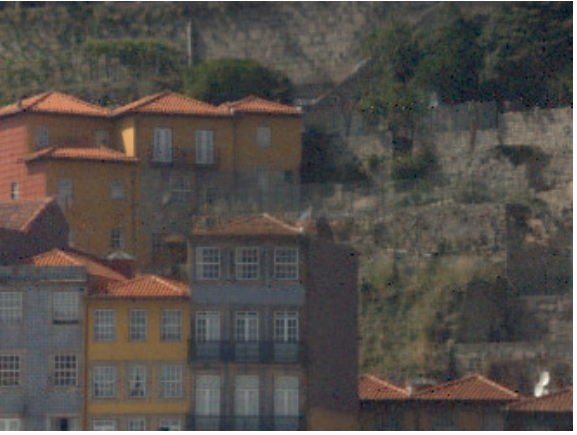}};
        \node[above] at ($(img)+(t)$) {\small TR-ALS};
        \node[below] at ($(img)-(d)$) {\scriptsize 38.3964dB};

        \node[inner sep=0] (img) at ($(O)+3*(C)$) {\includegraphics[width=1.69cm, height=1.25cm]{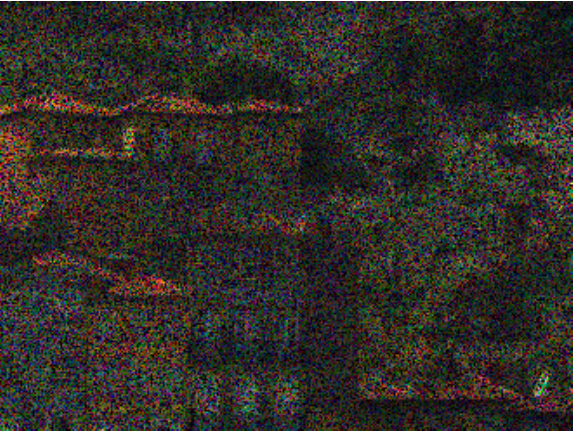}};
        \node[above] at ($(img)+(t)$) {\small HaLRTC};
        \node[below] at ($(img)-(d)$) {\scriptsize 21.3404dB};
        
        \node[inner sep=0] (img) at ($(O)+4*(C)$) {\includegraphics[width=1.69cm, height=1.25cm]{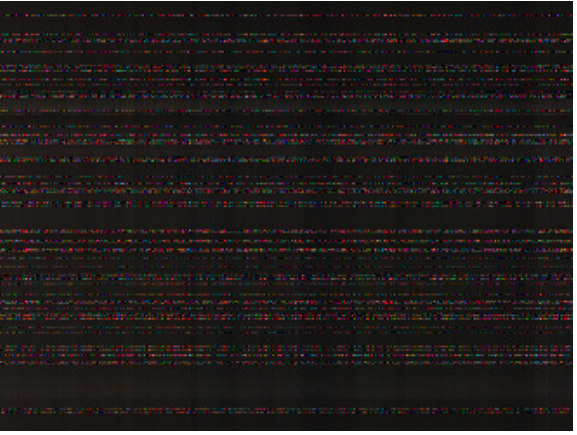}};
        \node[above] at ($(img)+(t)$) {\small TT-RCG};
        \node[below] at ($(img)-(d)$) {\scriptsize 19.3895dB};
        
        \node[inner sep=0] (img) at ($(O)+5*(C)+(0,0.003)$) {\includegraphics[width=1.69cm, height=1.24cm]{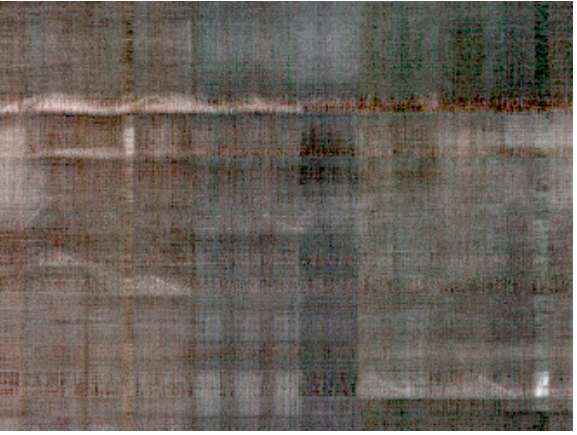}};
        \node[above] at ($(img)+(t)$) {\small CP-WOPT};
        \node[below] at ($(img)-(d)$) {\scriptsize 22.6050dB};
        
        \node[inner sep=0] (img) at ($(O)+6*(C)$) {\includegraphics[width=1.69cm, height=1.25cm]{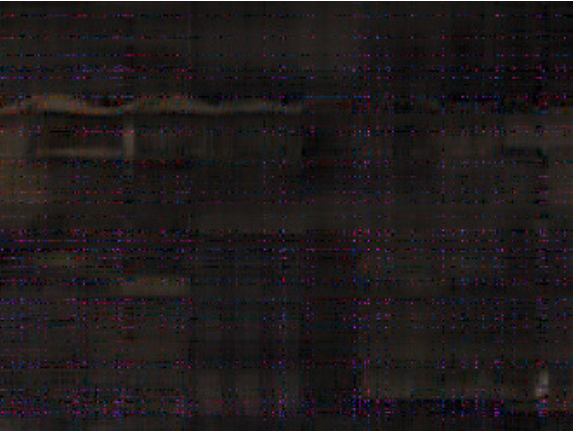}};
        \node[above] at ($(img)+(t)-(0,0.05)$) {\small geomCG};
        \node[below] at ($(img)-(d)$) {\scriptsize 19.9492dB};

        \coordinate (O) at (R);
        \node[inner sep=0] (image) at (O) {\includegraphics[width=1.69cm, height=1.25cm]{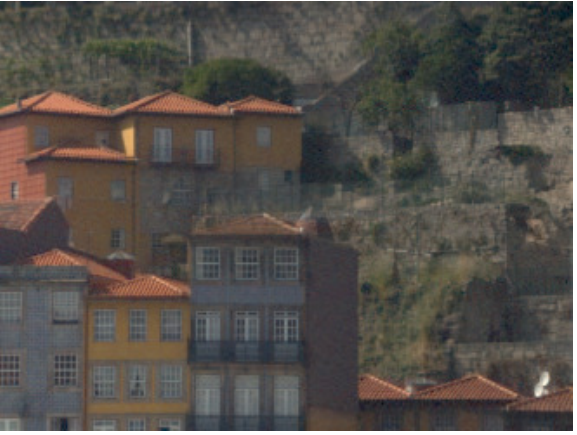}};
        \node[below] at ($(O)-(d)$) {\scriptsize 39.9493dB};
        
        \node[inner sep=0] (img) at ($(O)+(C)$) {\includegraphics[width=1.69cm, height=1.25cm]{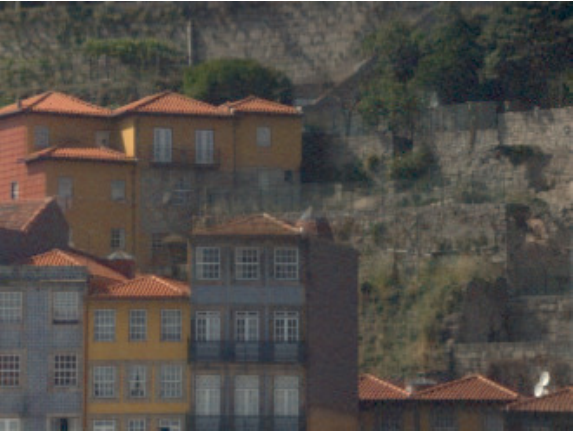}};
        \node[below] at ($(img)-(d)$) {\scriptsize 40.0437dB};
        
        \node[inner sep=0] (img) at ($(O)+2*(C)$) {\includegraphics[width=1.69cm, height=1.25cm]{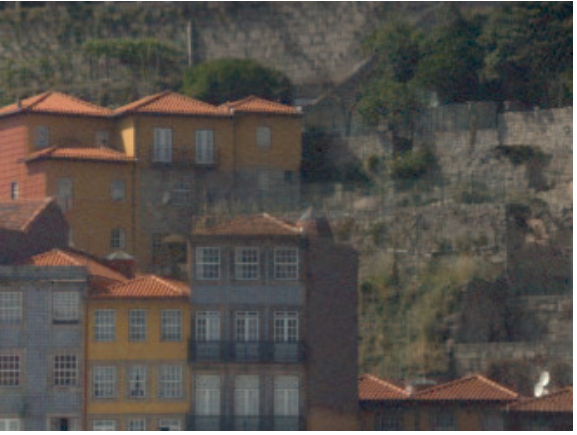}};
        \node[below] at ($(img)-(d)$) {\scriptsize 38.8803dB};

        \node[inner sep=0] (img) at ($(O)+3*(C)$) {\includegraphics[width=1.69cm, height=1.25cm]{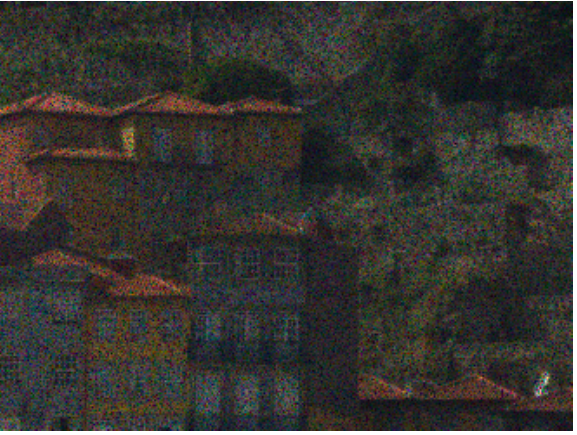}};
        \node[below] at ($(img)-(d)$) {\scriptsize 23.4366dB};
        
        \node[inner sep=0] (img) at ($(O)+4*(C)+(0,0.004)$) {\includegraphics[width=1.69cm, height=1.24cm]{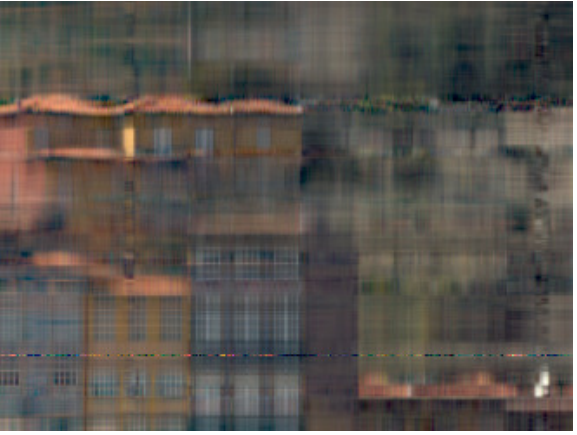}};
        \node[below] at ($(img)-(d)$) {\scriptsize 29.6771dB};
        
        \node[inner sep=0] (img) at ($(O)+5*(C)+(0,0.004)$) {\includegraphics[width=1.69cm, height=1.24cm]{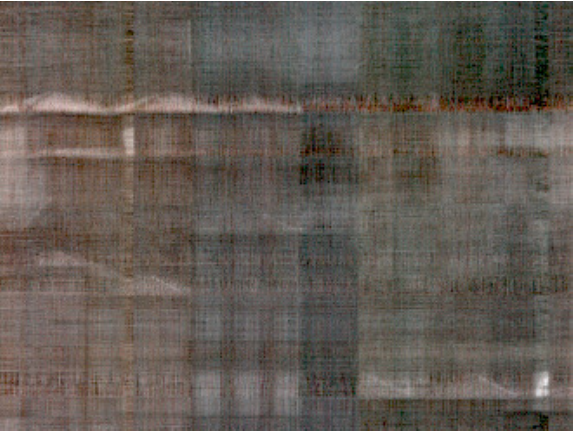}};
        \node[below] at ($(img)-(d)$) {\scriptsize 24.3917dB};
        
        \node[inner sep=0] (img) at ($(O)+6*(C)$) {\includegraphics[width=1.69cm, height=1.25cm]{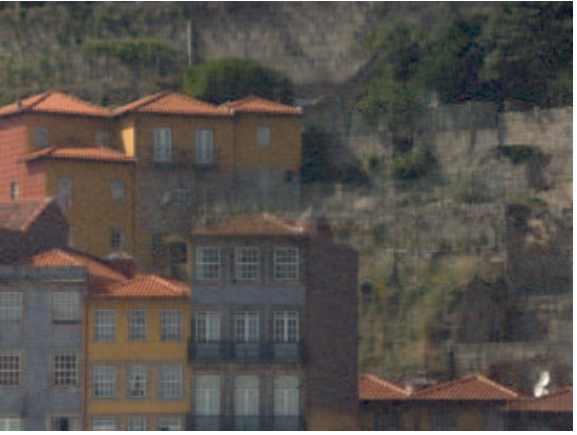}};
        \node[below] at ($(img)-(d)$) {\scriptsize 35.3738dB};

        \coordinate (O) at ($2*(R)$);
        \node[inner sep=0] (image) at (O) {\includegraphics[width=1.69cm, height=1.25cm]{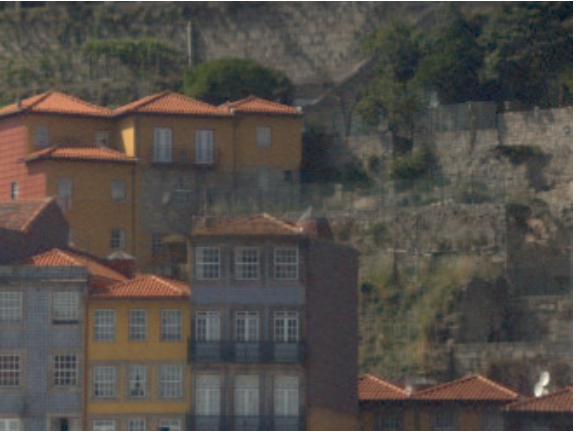}};
        \node[below] at ($(O)-(d)$) {\scriptsize 40.4171dB};
        
        \node[inner sep=0] (img) at ($(O)+(C)$) {\includegraphics[width=1.69cm, height=1.25cm]{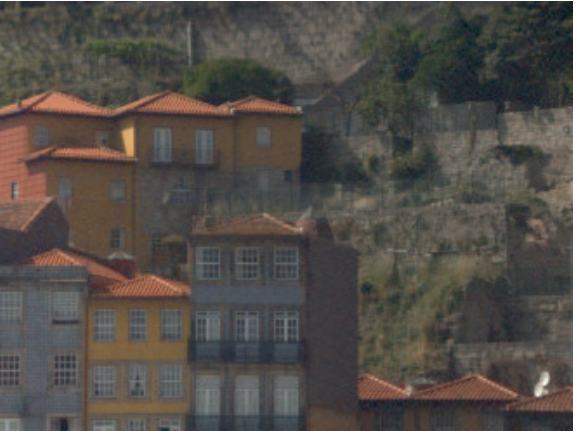}};
        \node[below] at ($(img)-(d)$) {\scriptsize 40.4871dB};
        
        \node[inner sep=0] (img) at ($(O)+2*(C)$) {\includegraphics[width=1.69cm, height=1.25cm]{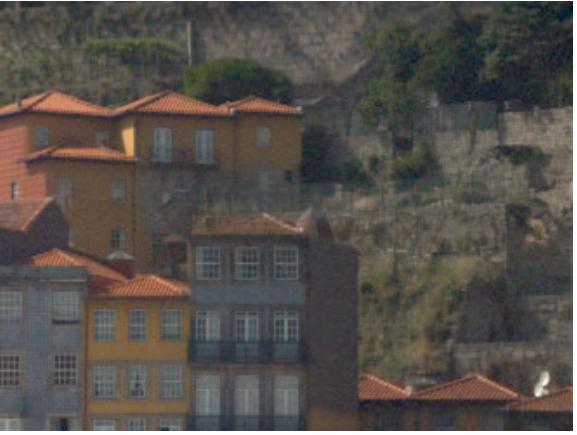}};
        \node[below] at ($(img)-(d)$) {\scriptsize 38.3704dB};

        \node[inner sep=0] (img) at ($(O)+3*(C)$) {\includegraphics[width=1.69cm, height=1.25cm]{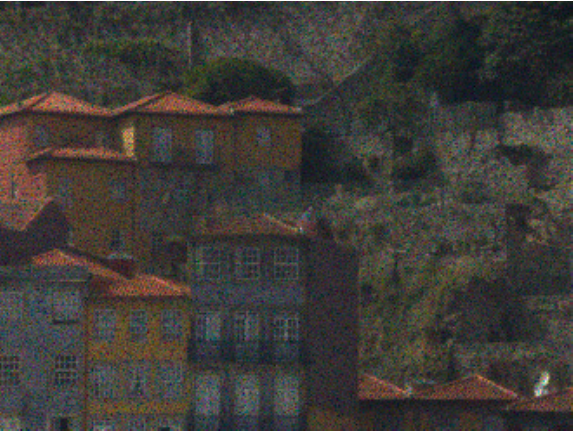}};
        \node[below] at ($(img)-(d)$) {\scriptsize 26.1920dB};
        
        \node[inner sep=0] (img) at ($(O)+4*(C)+(0,0.003)$) {\includegraphics[width=1.69cm, height=1.235cm]{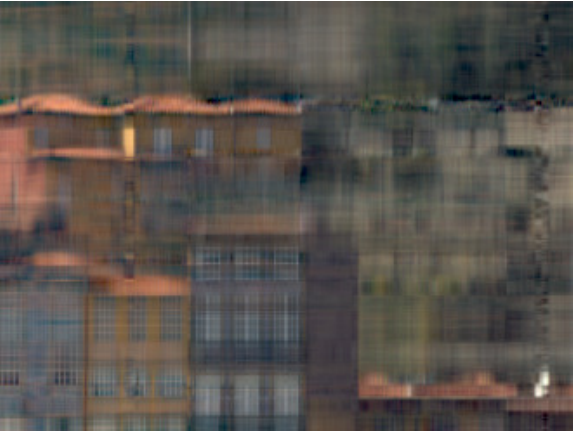}};
        \node[below] at ($(img)-(d)$) {\scriptsize 29.8119dB};
        
        \node[inner sep=0] (img) at ($(O)+5*(C)+(0,0.004)$) {\includegraphics[width=1.69cm, height=1.235cm]{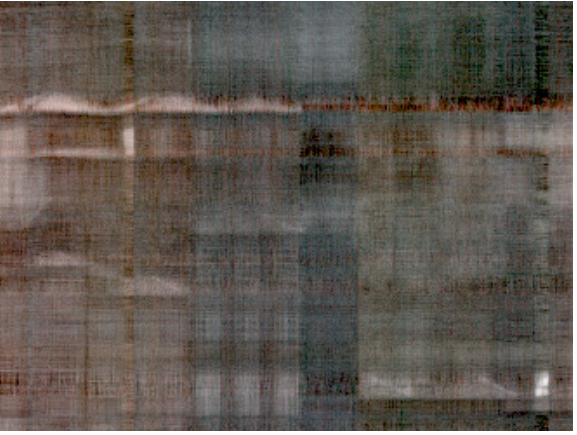}};
        \node[below] at ($(img)-(d)$) {\scriptsize 25.7736dB};
        
        \node[inner sep=0] (img) at ($(O)+6*(C)$) {\includegraphics[width=1.69cm, height=1.25cm]{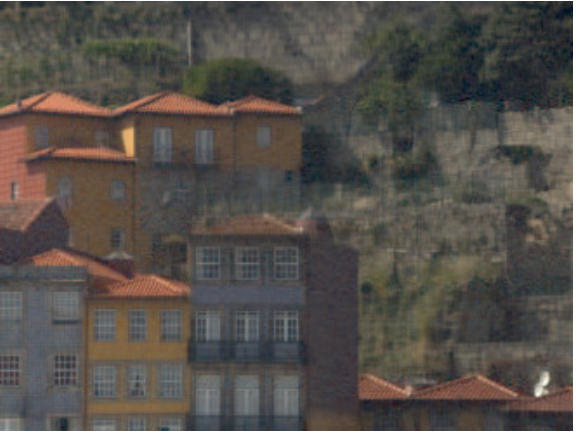}};
        \node[below] at ($(img)-(d)$) {\scriptsize 36.2973dB};

        \coordinate (O) at ($3*(R)$);
        \node[inner sep=0] (image) at (O) {\includegraphics[width=1.69cm, height=1.25cm]{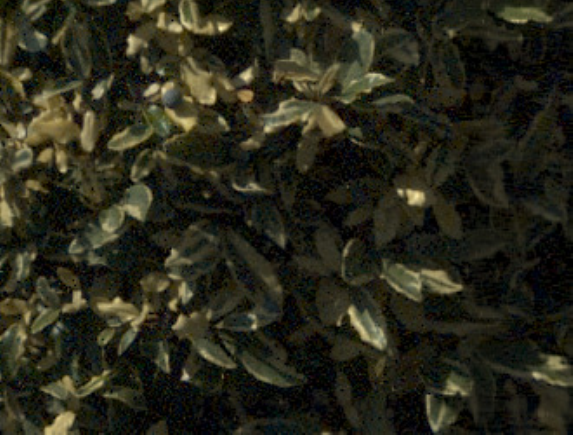}};
        \node[below] at ($(O)-(d)$) {\scriptsize 40.5529dB};
        
        \node[inner sep=0] (img) at ($(O)+(C)$) {\includegraphics[width=1.69cm, height=1.25cm]{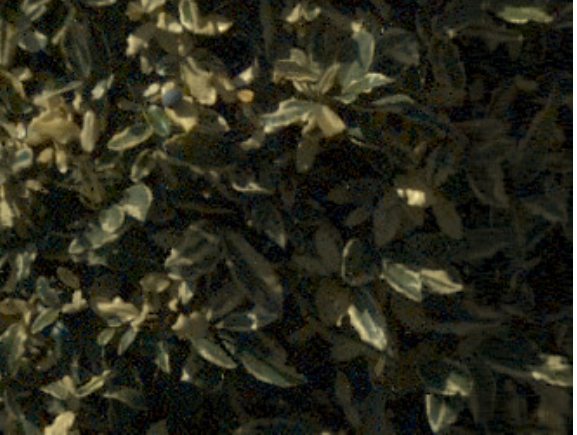}};
        \node[below] at ($(img)-(d)$) {\scriptsize 40.0065dB};
        
        \node[inner sep=0] (img) at ($(O)+2*(C)$) {\includegraphics[width=1.69cm, height=1.25cm]{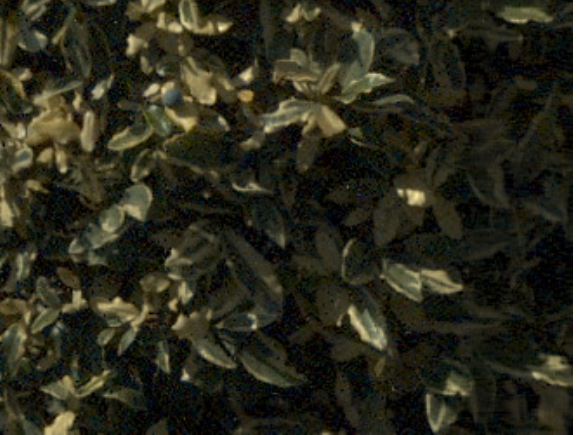}};
        \node[below] at ($(img)-(d)$) {\scriptsize 40.1682dB};

        \node[inner sep=0] (img) at ($(O)+3*(C)$) {\includegraphics[width=1.69cm, height=1.25cm]{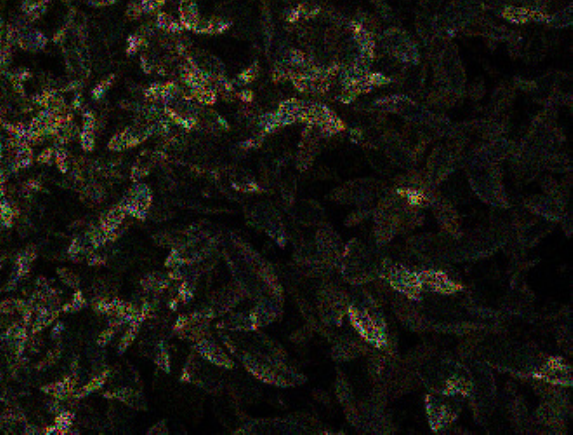}};
        \node[below] at ($(img)-(d)$) {\scriptsize 24.7499dB};
        
        \node[inner sep=0] (img) at ($(O)+4*(C)$) {\includegraphics[width=1.69cm, height=1.25cm]{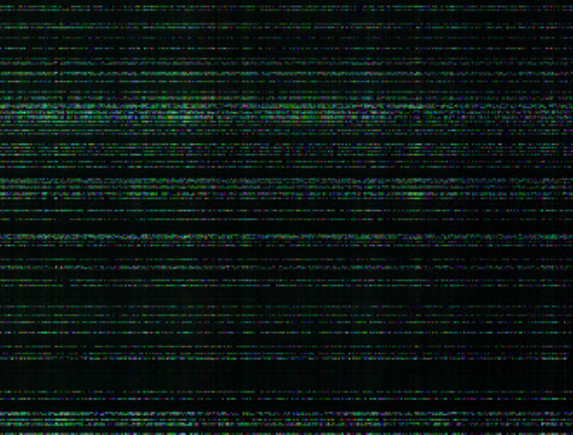}};
        \node[below] at ($(img)-(d)$) {\scriptsize 21.4320dB};
        
        \node[inner sep=0] (img) at ($(O)+5*(C)$) {\includegraphics[width=1.69cm, height=1.252cm]{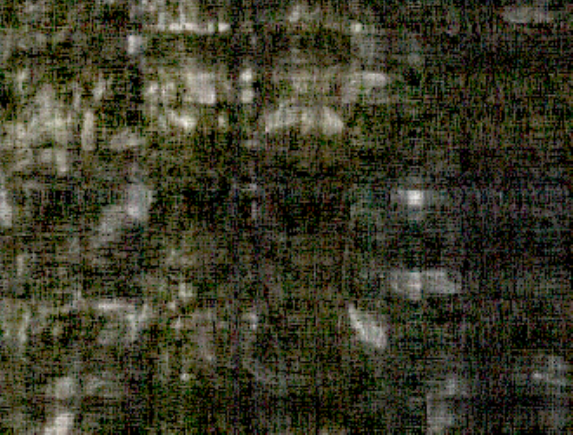}};
        \node[below] at ($(img)-(d)$) {\scriptsize 23.5259dB};
        
        \node[inner sep=0] (img) at ($(O)+6*(C)$) {\includegraphics[width=1.69cm, height=1.25cm]{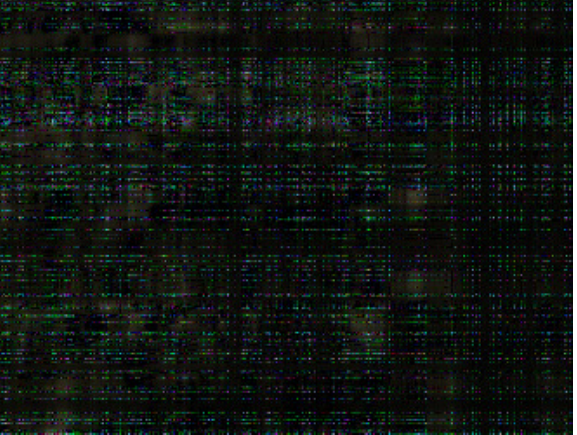}};
        \node[below] at ($(img)-(d)$) {\scriptsize 21.9541dB};

        \coordinate (O) at ($4*(R)$);
        \node[inner sep=0] (image) at (O) {\includegraphics[width=1.69cm, height=1.25cm]{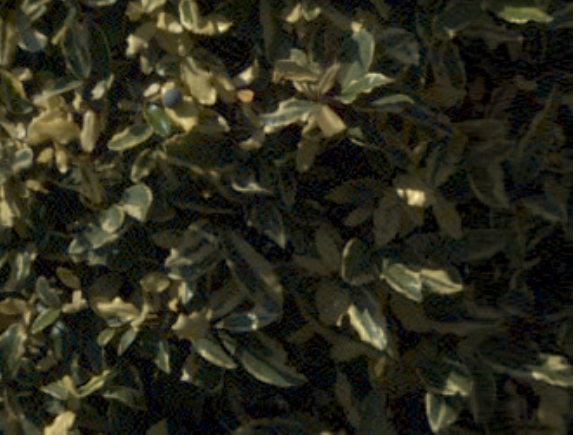}};
        \node[below] at ($(O)-(d)$) {\scriptsize 40.4514dB};
        
        \node[inner sep=0] (img) at ($(O)+(C)$) {\includegraphics[width=1.69cm, height=1.25cm]{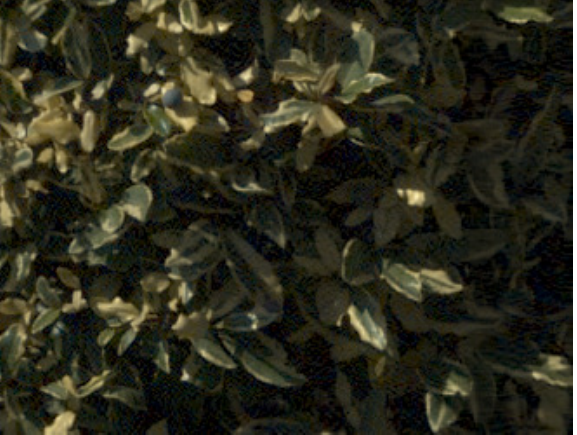}};
        \node[below] at ($(img)-(d)$) {\scriptsize 40.3119dB};
        
        \node[inner sep=0] (img) at ($(O)+2*(C)$) {\includegraphics[width=1.69cm, height=1.25cm]{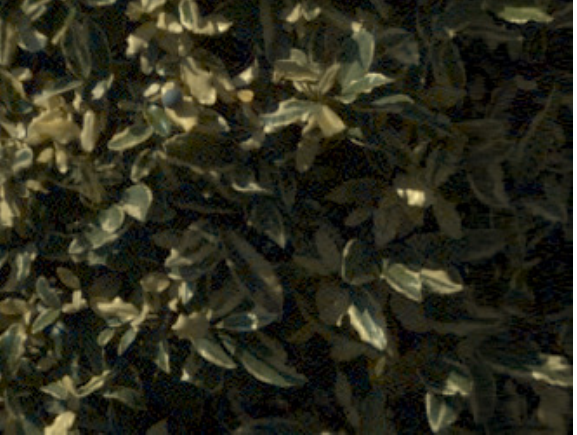}};
        \node[below] at ($(img)-(d)$) {\scriptsize 39.9008dB};

        \node[inner sep=0] (img) at ($(O)+3*(C)$) {\includegraphics[width=1.69cm, height=1.25cm]{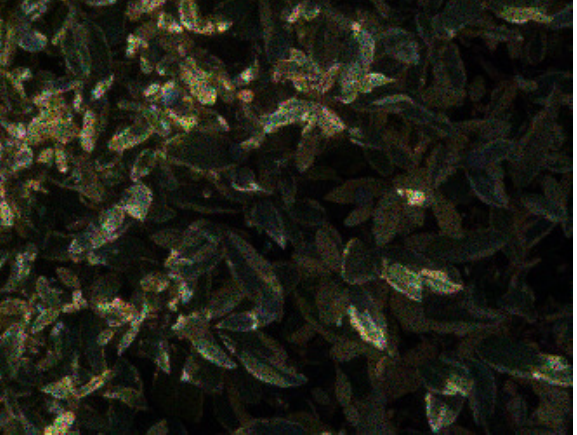}};
        \node[below] at ($(img)-(d)$) {\scriptsize 27.6014dB};
        
        \node[inner sep=0] (img) at ($(O)+4*(C)$) {\includegraphics[width=1.69cm, height=1.25cm]{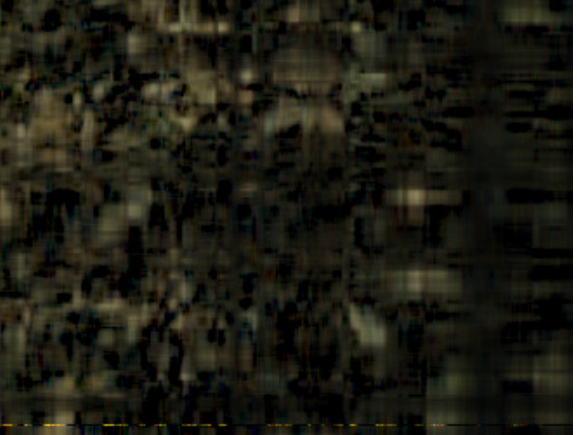}};
        \node[below] at ($(img)-(d)$) {\scriptsize 25.3039dB};
        
        \node[inner sep=0] (img) at ($(O)+5*(C)$) {\includegraphics[width=1.69cm, height=1.252cm]{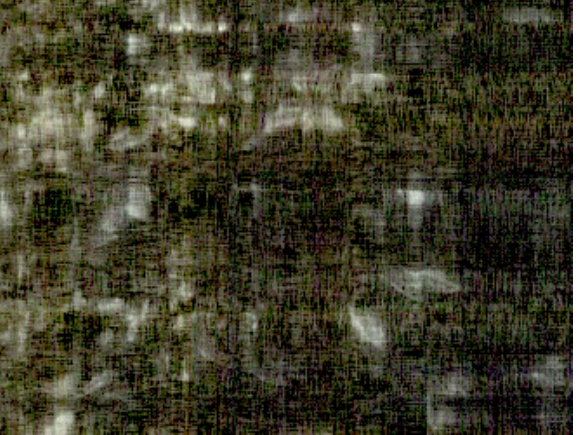}};
        \node[below] at ($(img)-(d)$) {\scriptsize 23.1848dB};
        
        \node[inner sep=0] (img) at ($(O)+6*(C)$) {\includegraphics[width=1.69cm, height=1.25cm]{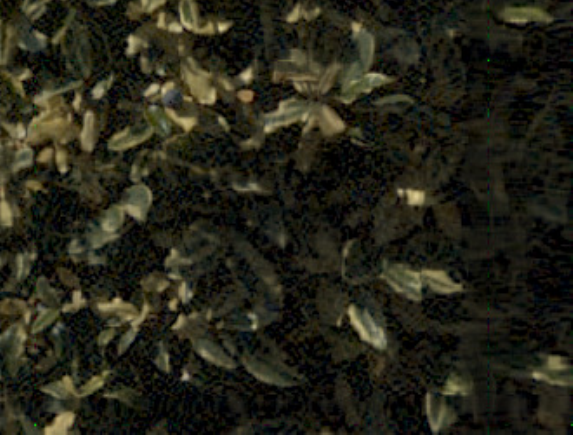}};
        \node[below] at ($(img)-(d)$) {\scriptsize 35.0618dB};

        \coordinate (O) at ($5*(R)$);
        \node[inner sep=0] (image) at (O) {\includegraphics[width=1.69cm, height=1.25cm]{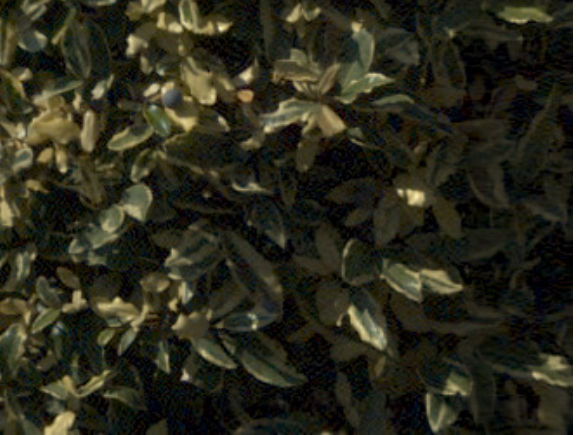}};
        \node[below] at ($(O)-(d)$) {\scriptsize 40.4516dB};
        
        \node[inner sep=0] (img) at ($(O)+(C)$) {\includegraphics[width=1.69cm, height=1.25cm]{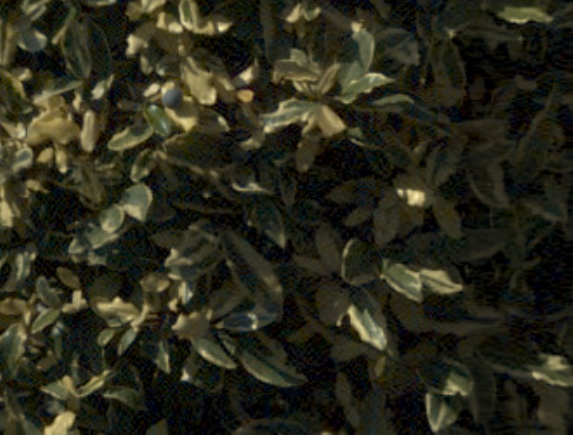}};
        \node[below] at ($(img)-(d)$) {\scriptsize 39.9737dB};
        
        \node[inner sep=0] (img) at ($(O)+2*(C)$) {\includegraphics[width=1.69cm, height=1.25cm]{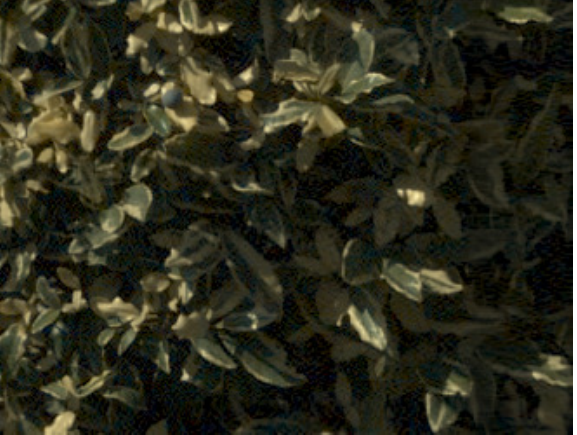}};
        \node[below] at ($(img)-(d)$) {\scriptsize 39.4380dB};

        \node[inner sep=0] (img) at ($(O)+3*(C)$) {\includegraphics[width=1.69cm, height=1.25cm]{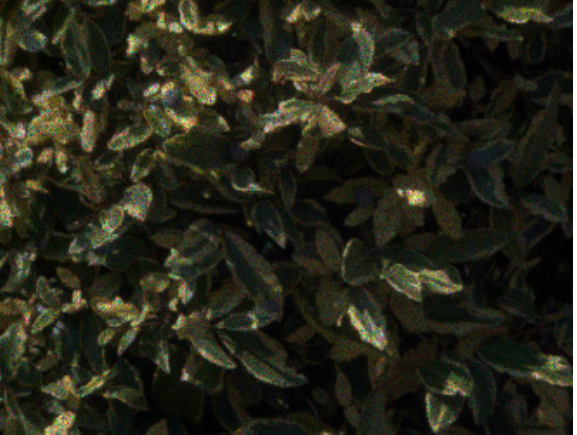}};
        \node[below] at ($(img)-(d)$) {\scriptsize 30.1986dB};
        
        \node[inner sep=0] (img) at ($(O)+4*(C)+(-0.0035,0)$) {\includegraphics[width=1.682cm, height=1.25cm]{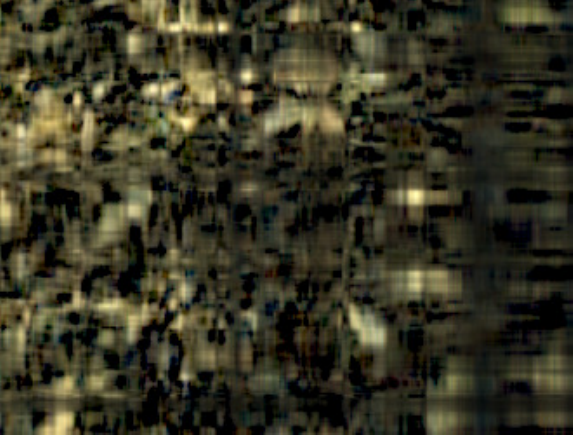}};
        \node[below] at ($(img)-(d)$) {\scriptsize 24.9847dB};
        
        \node[inner sep=0] (img) at ($(O)+5*(C)+(0.0045,0)$) {\includegraphics[width=1.70cm, height=1.25cm]{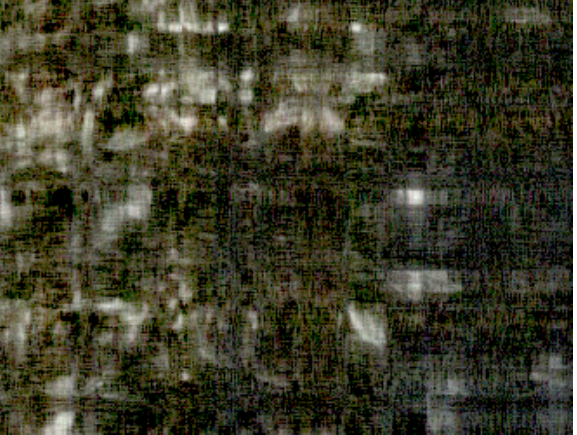}};
        \node[below] at ($(img)-(d)$) {\scriptsize 23.2231dB};
        
        \node[inner sep=0] (img) at ($(O)+6*(C)$) {\includegraphics[width=1.69cm, height=1.25cm]{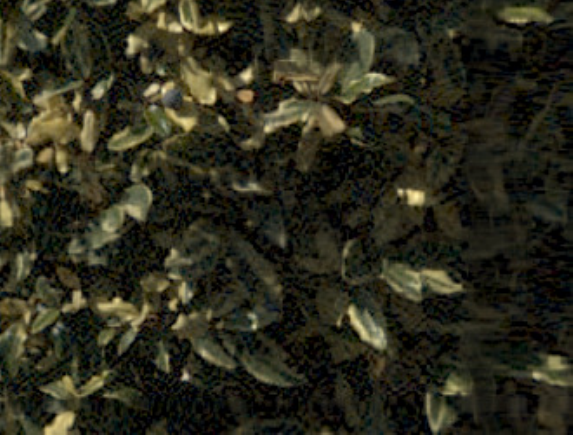}};
        \node[below] at ($(img)-(d)$) {\scriptsize 34.9976dB};
    \end{tikzpicture}

    \caption{RGB representations of recovered images by different completion algorithms. The first three rows represent recovery results of the ``Ribeira'' image, under sampling rates $p=0.1,0.3,0.5$ in each row. The last three rows represent recovery results of the ``Bush'' image, under sampling rates $p=0.1,0.3,0.5$ in each row. The PSNR is displayed under each image 
    }
    \label{fig: results for HSIs}
\end{figure}

Numerical results for the completion of two hyperspectral images are shown in Fig.~\ref{fig: results for HSIs} and Table~\ref{tab: results for HSIs}. One difficulty of recovering the ``Ribeira'' image is shrubs around the buildings, which triggers different recovery quality of different algorithms. The recovery quality of ``Bush'' image depends on the details of leaves. Figure~\ref{fig: results for HSIs} displays the recovery results under different sampling rates $p=0.1,0.3,0.5$. The images recovered by TR-based algorithms (TR-RGD, TR-RCG, TR-ALS) depict more details than other candidates, e.g., the shrubs and the bushes. In fact, this observation can be quantified by the PSNR and relative errors; see Table~\ref{tab: results for HSIs}. The proposed methods are favorably comparable to TR-ALS. Specifically, we observe that the proposed TR-RCG algorithm reaches the highest PSNR and the lowest relative error among all compared methods in most cases.

\begin{table}[htbp]
    \centering
    \caption{PSNR and relative errors for completion of two hyperspectral images}
    \scriptsize
    \resizebox{\textwidth}{!}{
    \begin{tabular}{ccrrrrrrr}
        \toprule
        \multirow{2}*{$p$} & \multirow{2}*{Results} &TR-RCG & TR-RGD & TR-ALS & HaLRTC & TT-RCG & CP-WOPT & geomCG\\
        \cmidrule{3-9}
        & &  \multicolumn{7}{c}{Ribeira House Shrubs}\\
        \midrule
        \multirow{2}*{$0.1$} & PSNR & 38.4310 & 37.5599 & 38.3964 & 21.3404 & 19.3895 & 22.6050 & 19.9492 \\ 
        & relerr & 0.1058 & 0.1170 & 0.1062 & 0.7569 & 0.9476 & 0.6544 & 0.8884 \\ 
        % \midrule
        \multirow{2}*{$0.3$}  & PSNR & 39.9493 & 40.0437 & 38.8803 & 23.4366 & 29.6771 & 24.3917 & 35.3738 \\ 
         & relerr & 0.0888 & 0.0879 & 0.1005 & 0.5946 & 0.2899 & 0.5327 & 0.1505 \\ 
         \multirow{2}*{$0.5$}  & PSNR & 40.4171 & 40.4871 & 38.3704 & 26.1920 & 29.8119 & 25.7736 & 36.2973 \\ 
         & relerr & 0.0842 & 0.0835 & 0.1066 & 0.4330 & 0.2854 & 0.4544 & 0.1353 \\ 
        \midrule
         &  & \multicolumn{7}{c}{Bom Jesus Bush}\\
        % \cmidrule{3-9}
        % & & TR-RCG & TR-RGD & TR-ALS & HaLRTC & TT-RCG & CP-WOPT & geomCG\\
        \midrule
        \multirow{2}*{$0.1$} & PSNR & 40.5529 & 40.0065 & 40.1682 & 24.7499 & 21.4320 & 23.5259 & 21.9541 \\ 
         & relerr & 0.1185 & 0.1262 & 0.1239 & 0.7310 & 1.0711 & 0.8417 & 1.0086 \\ 
         \multirow{2}*{$0.3$} & PSNR & 40.4514 & 40.3119 & 39.9008 & 27.6014 & 25.3039 & 23.1848 & 35.0618 \\ 
         & relerr & 0.1199 & 0.1219 & 0.1278 & 0.5265 & 0.6859 & 0.8754 & 0.2230 \\ 
        %  \midrule
         \multirow{2}*{$0.5$} & PSNR & 40.4516 & 39.9737 & 39.4380 & 30.1986 & 24.9847 & 23.2231 & 34.9976 \\ 
         & relerr & 0.1199 & 0.1267 & 0.1348 & 0.3904 & 0.7115 & 0.8715 & 0.2247 \\ 
        \bottomrule
    \end{tabular}
    }
    \label{tab: results for HSIs}
\end{table}

\subsection{Experiments on high-dimensional functions}
The discretization of high-dimensional functions ${h:[0,1]^d\to\mathbb{R}}$ requires storing a large tensor $\tensA\in\mathbb{R}^{n_1\times n_2\times\cdots\times n_d}$, which is unfavorable in practice. To this end, one can apply tensor decomposition and tensor completion to recover the tensor~$\tensA$ by partially observed entries; see applications in~\cite{glau2020low}. In this subsection, we compare TR-based algorithms (TR-RGD, TR-RCG, TR-ALS) with TT-RCG on completion of the data tensor $\tensA$ generated from $h$ by evenly dividing $[0,1]$ in dimension $k$ of~$[0,1]^d$ into $n_k-1$ intervals for $k\in[d]$. Specifically, \[\tensA(i_1,i_2,\dots,i_d)=h\left(\frac{i_1-1}{n_1-1},\frac{i_2-1}{n_2-1},\dots,\frac{i_d-1}{n_d-1}\right),\ i_k\in[n_k],\ k\in[d].\] 
We consider the following two functions~\cite[Sect. 5.4]{steinlechner2016riemannian},
\begin{equation*}
    \begin{aligned}
        h_1\!:&\ [0,1]^d\to\mathbb{R},\quad h_1(\mathbf{x}){:=}\exp(-\|\mathbf{x}\|),\quad \text{and}\\
        h_2\!:&\ [0,1]^d\to\mathbb{R},\quad h_2(\mathbf{x}){:=}\frac{1}{\|\mathbf{x}\|},
    \end{aligned}
\end{equation*}
We set $d=4,\ n_1=n_2=n_3=n_4=20$, and the sampling rate $p=0.001, 0.005,$ $ 0.01, 0.05, 0.1$. We follow the same way in subsection~\ref{subsec: synthetic} to create the sampling set~$\Omega$. The rank-increasing strategy is implemented to both TT and TR algorithms by following~\cite{steinlechner2016riemannian}. In order to choose a search space with similar size, we set the maximum rank of TT-RCG as $(1,5,5,5,1)$, and the maximum rank for TR-based algorithms as $(4,4,4,4)$; see Appendix~\ref{app: search space} for details. We adopt the stopping criteria in section~\ref{sec: numerical exps}. Additionally, due to the rank-increasing strategy, an algorithm is also terminated if: 1) maximum iteration number 50 is reached in each fixed-rank searching; 2) the maximum rank is achieved; 3) there is no acceptance of rank increase along any mode for a given point.

\begin{table}[htbp]
	\centering
	\caption{Test errors for high-dimensional functions}
	\footnotesize
    \resizebox{\textwidth}{!}{
	\begin{tabular}{lcccccccc}
		\toprule
		\multicolumn{1}{c}{\multirow{2}*{$p$}}& \multicolumn{4}{c}{$\exp(-\|\mathbf{x}\|)$} & \multicolumn{4}{c}{${1}/{\|\mathbf{x}\|}$}\\
		\cmidrule(lr){2-5} \cmidrule(lr){6-9}
		& TR-RGD & TR-RCG & TR-ALS & TT-RCG & TR-RGD & TR-RCG & TR-ALS & TT-RCG\\
		\midrule
		{0.001} & 8.0884e-2 & 7.4157e-2 & 7.4161e-2 & 1.3445e-1 & 1.7531e-1 & 1.8106e-1 & 1.8081e-1 & 2.6876e-1\\
		{0.005} & 7.3505e-3 & 8.7366e-3 & 9.2121e-3 & 1.5904e-2 & 3.4428e-2 & 2.9218e-2 & 3.2090e-2 & 1.2899e-1\\
		{0.01} & 6.2650e-3 & 9.7247e-4 & 1.8737e-3 & 4.1233e-3 & 2.5230e-2 & 1.7676e-2 & 1.8697e-2 & 3.4675e-2\\
		{0.05} & 3.8862e-4 & 1.5019e-4 & 1.8218e-4 & 2.2991e-4 & 3.8510e-3 & 3.6002e-3 & 5.2173e-3 & 3.5697e-3\\
		{0.1} & 1.2251e-4 & 5.9871e-5 & 6.8898e-5 & 8.2512e-5 & 7.7886e-4 & 2.9423e-4 & 6.0423e-4 & 7.4727e-4\\
		\bottomrule
	\end{tabular}
    }
	\label{tab: Interp of funcs}
\end{table}

The numerical results for recovering high-dimensional functions $h_1$ and $h_2$ are reported in Table~\ref{tab: Interp of funcs}. It illustrates that all algorithms have comparable performance. TR-based algorithms perform favorably comparable to TT-RCG. Among all TR-based algorithms, TR-RCG has a better performance in most experiments.

\section{Conclusion and perspectives}\label{sec: conclusion}
We have developed Riemannian preconditioned algorithms for the tensor completion problem based on tensor ring decomposition. The preconditioning effect stems from a metric defined on the product space of matrices generated from the mode-2 unfolding of core tensors. However, the straightforward calculation of Riemannian gradient requires large matrix multiplications that are unaffordable in practice. To this end, we have adopted a procedure that efficiently computes the Riemannian gradient without forming large matrices explicitly. The proposed algorithms enjoy global convergence results. Numerical comparisons on both synthetic and real-world datasets present promising results. 

Since the rank parameter has to be fixed a priori, we are interested in rank-adaptive strategies that may result in more accurate completion.

% \begin{appendices}
\appendix
     \section{Speedup of efficient gradient computation} \label{app: speedup} In contrast with naive computation, Algorithm~\ref{alg: compute uneqk t uneqk} provides an efficient way to compute the gradients. Table~\ref{tab: Speedups} reports the speedup results on MovieLens 1M dataset. ``Naive'' denotes the non-optimized implementation, i.e., explicitly forming $\matW_{\neq k}$ and large matrix multiplication. The results show that the proposed algorithm is of great potential to be applied in large-scale problems. 

    \begin{table}[htbp]
        \centering
        \caption{Speedup of efficient gradient computation on MovieLens 1M dataset. ``Avg. Speedup'': average speedup}
        \scriptsize
        \begin{tabular}{rrrccccc}
            \toprule
            \multirow{2}*{\#iter} & \multicolumn{2}{c}{Time (seconds)} & \multirow{2}*{Avg. Speedup} &\multicolumn{2}{c}{Relerr on $\Omega$ ($\varepsilon_{\Omega}$)} & \multicolumn{2}{c}{Relerr on $\Gamma$ ($\varepsilon_{\Gamma}$)}\\
            \cmidrule(lr){2-3} \cmidrule(lr){5-6} \cmidrule(lr){7-8}
            & \multicolumn{1}{c}{Naive}  & \multicolumn{1}{c}{Proposed} &  & Naive & Proposed & Naive & Proposed\\
            \midrule
            1 & 266.6065 & 5.9717 & - & 4.0151 & 4.0151 & 4.0233 & 4.0233\\
            21 & 4880.3810 & 37.7889 & 145.5539$\times$ & 0.7636 & 0.7636 & 0.8674 & 0.8674\\
            41 & 9400.1950 & 69.2155 & 144.6929$\times$ & 0.3979 & 0.3979 & 0.5268 & 0.5268\\
            61 & 13957.4616 & 100.5884 & 144.8814$\times$ & 0.2796 & 0.2796 & 0.4084 & 0.4084\\
            81 & 18518.4900 & 132.6191 & 144.2527$\times$ & 0.2612 & 0.2612 & 0.3880 & 0.3880\\
            \bottomrule
        \end{tabular}
        \label{tab: Speedups}
    \end{table}

        \section{The rank selection for numerical experiments}\label{app: search space}

        Table~\ref{tab: rank selections} introduces the rank selections in numerical experiments. For instance, we initially choose the Tucker rank $(65,65,7)$ following~\cite{kressner2014low} in experiments on hyperspectral images with size $250\times 330\times 33$. In order to ensure a fair comparison, the search spaces for different tensor formats are of a similar size. To this end, we compute the number of parameters by
        \[65\times 65\times 7+250\times 65+330\times 65+33\times 7=67506.\]
        For TT decomposition, we select the rank as $(1,15,15,1)$ with 
        \[250\times 15+15\times 330\times 15+15\times 33=78495\]
        parameters. For TR decomposition, the rank is selected by $(7,16,7)$. Then, the number of parameters is 
        \[7\times 250\times 16+16\times 330\times 7+7\times 33\times 7=66577.\]
        For CP decomposition, the rank is selected by
        \[R=\frac{67506}{250+ 330+ 33}\approx 110.\]
        Therefore, the search spaces appear to have comparable sizes across various tensor formats. The ranks are chosen using the same procedure in other experiments.

        \begin{table}[htbp]
            \centering
            \caption{Rank selections for different tensor formats in numerical experiments. ``\#params'': number of parameters}
            \label{tab: rank selections}

            \resizebox{\textwidth}{!}{
            \begin{tabular}{ccrcrcrcr}
                \toprule
                \multirow{2}*{Format} & \multicolumn{4}{c}{MovieLens 1M} & \multicolumn{2}{c}{Hyperspectral images} & \multicolumn{2}{c}{High-dimensional functions}\\
                \cmidrule(lr){2-3} \cmidrule(lr){4-5} \cmidrule(lr){6-7} \cmidrule(lr){8-9}
                 & Rank & \#params & Rank & \#params & Rank & \#params & Rank & \#params\\
                \midrule
                CP & 36 & 365112 & 49 & 496958 & 110 & 67430 & - & - \\
                Tucker & (36,36,36) & 411768 & (60,30,18) & 516060 & (65,65,7) & 67506 & - & - \\
                TT & (1,9,9,1) & 375822 & (1,10,10,1) & 457100 & (1,15,15,1) & 78495 & (1,5,5,5,1) & 1200\\
                TR & (6,6,6) & 365112 & (6,10,3) & 483660 & (7,16,7) & 66577 & (4,4,4,4) & 1280\\
                \bottomrule
            \end{tabular}
            }
        \end{table}




% \end{appendices}

    \section*{Acknowledgements}
We would like to thank the two anonymous reviewers for helpful comments. We acknowledge Shuyu Dong for helpful discussions on the preconditioned metric.

\section*{Declaration}
The authors declare that the data supporting the findings of this study are available within the paper. The authors have no competing interests to declare that are relevant to the content of this article.

%\begin{acknowledgements}
%\end{acknowledgements}

\bibliographystyle{spmpsci}      % mathematics and physical sciences
\bibliography{sn-bibliography}   % name your BibTeX data base

% \begin{thebibliography}{}
% % Format for Journal Reference
% Author, Article title, Journal, Volume, page numbers (year)
% Author, Book title, page numbers. Publisher, place (year)
% \end{thebibliography}

\end{document}